\documentclass[12pt]{article}
\usepackage{latexsym,amsmath,amssymb}
\usepackage[margin=1.2in]{geometry}

\newtheorem{thm}{Theorem}[section]

\newtheorem{lemma}[thm]{Lemma}

\newtheorem{conj}[thm]{Conjecture}

\newtheorem{claim}[thm]{Claim}

\newcommand{\ul}[0]{\underline}
\newcommand{\qed}[0]{{\hspace*{\fill}\mbox{$\Box$}}}
\newcommand{\qq}{[q]}

\begin{document}
\renewcommand{\thefootnote}{\fnsymbol{footnote}}

\title{Phase coexistence and torpid mixing in the $3$-coloring model on ${\mathbb Z}^d$}

\author{David Galvin\thanks{Department of Mathematics, University of Notre Dame,
Notre Dame IN, USA; dgalvin1@nd.edu.},~~Jeff Kahn\thanks{Department of Mathematics, Rutgers University, New Brunswick NJ, USA; jkahn@math.rutgers.edu.},~~Dana
Randall\thanks{School of Computer Science, Georgia Institute of Technology,
Atlanta GA, USA; randall@cc.gatech.edu.},~~Gregory B.
Sorkin\thanks{Departments of Management and Mathematics, London School of Economics, London, England; g.b.sorkin@lse.ac.uk.}}

\maketitle

% DJG October 10 2012: changed slow to torpid throughout

\begin{abstract}
We show that for all sufficiently large $d$, the uniform proper $3$-coloring
model (in physics called the $3$-state antiferromagnetic Potts model at zero temperature) on ${\mathbb Z}^d$ admits multiple maximal-entropy Gibbs measures. This is a consequence of the following combinatorial result: if a proper $3$-coloring is chosen uniformly from a box in ${\mathbb Z}^d$, conditioned on color $0$ being given to all
the vertices on the boundary of the box which are at an odd distance from a
fixed vertex $v$ in the box, then the probability that $v$
gets color $0$ is exponentially small in $d$.

The proof proceeds through an analysis of a certain type of cutset separating $v$ from the boundary of the box, and builds on techniques developed by Galvin and Kahn in their proof of phase transition in the hard-core model on ${\mathbb Z}^d$.

Building further on these techniques, we study local Markov chains for sampling proper $3$-colorings of the discrete torus
${\mathbb Z}^d_n$.
We show that there is a constant $\rho \approx 0.22$ such that for all even $n \geq 4$ and $d$ sufficiently large, if ${\mathcal M}$ is a Markov chain on the set of proper $3$-colorings
of ${\mathbb Z}^d_n$ that updates the color of at most $\rho  n^d$ vertices
at each step and whose stationary distribution is uniform, then the
mixing time of ${\mathcal M}$ (the time taken for ${\mathcal M}$ to reach a distribution that is close to uniform, starting from an arbitrary coloring) is essentially exponential in $n^{d-1}$.
%\footnote{The mixing result was presented by the first and third authors at the 2007 ACM-SIAM Symposium on Discrete Algorithms (SODA '07) \cite{GalvinRandall}}.
\end{abstract}

\section{Introduction} \label{sec-intro}

A {\em(proper) $q$-coloring} of a graph $G=(V,E)$ is a function
$\chi:V(G)\rightarrow \qq$ satisfying $\chi(u) \neq \chi(v)$ whenever $uv \in
E$, where we use the notation $\qq=\{0, \ldots, q-1\}$. In the language of
statistical physic, a $q$-coloring of $G$ is a configuration in the {\em
zero-temperature $q$-state antiferromagnetic Potts model} on $G$ \cite{Potts}.
This is a simple model of the occupation of space by a collection of $q$
types of particles: % self-repulsive particles, with
the vertices of $G$ represent sites, each
occupied by exactly one particle, and the edges of $G$ represent pairs of sites
that are bonded (by spatial proximity, for example) and cannot be occupied by
particles of the same type. We write ${\mathcal C}_q(G)$, or simply ${\mathcal C}_q$, for the set of $q$-colorings of $G$.

%\gbs{Why not ``a'' or ``the'' uniform measure? Usage here eschews articles,
%consistently, and also for ``phase transition''.}
A basic question concerning ${\mathcal C}_q$ is, what does a typical (uniformly chosen) element look like? For finite $G$, uniform measure on ${\mathcal C}_q$ is unambiguous. For infinite $G$, the standard approach to defining uniform measure on ${\mathcal C}_q$ is through the notion of a {\em Gibbs measure}, which, roughly speaking, is a measure on ${\mathcal C}_q$ whose restriction to any finite subset of $V$ is uniform.

Formally, for finite $W \subseteq V$ let $\mu_{W^+}$ be uniform measure on the subgraph of $G$ induced by $W \cup \partial_{\rm ext} W$, where $\partial_{\rm ext} W$ is the set of vertices outside $W$ that are adjacent to something in $W$. Let $\mu$ be a measure on $({\mathcal C}_q,{\mathcal F}_{\rm cyl})$, where ${\mathcal F}_{\rm cyl}$ is the
$\sigma$-algebra generated by the cylinder events $\{\chi(v)=i\}$
for $v\in V$ and $i\in \qq$. We say that $\mu$ is a Gibbs measure {\em (with
uniform specification)} if the following condition holds: for all finite $W
\subseteq V$ and $\mu$-almost-all $\chi \in {\mathcal C}_q$, the probability
that $\chi'$ agrees with $\chi$ on $W$ given that it agrees with $\chi$ off
$W$, with $\chi'$ drawn according to $\mu$, is the same as the probability
that $\chi'$ agrees with $\chi$ on $W$ given that it agrees with $\chi$ on
$\partial_{\rm ext} W$, with $\chi'$ drawn according $\mu_{W^+}$. (See e.g.\ \cite{Georgii2} for a thorough treatment of this topic).

General compactness arguments show that an infinite graph $G$ admits at least one Gibbs measure.
A simple recipe for producing one is the following. For $\chi \in {\mathcal C}_q={\mathcal C}_q(G)$ and $W \subseteq V$, let ${\mathcal C}_q^\chi(W)$ be the set of colorings that agree with $\chi$ off $W$.
Fix $\chi \in {\mathcal C}_q$ and a nested sequence $(W_i)_{i=1}^\infty$ of
finite subsets of $V$ satisfying $\cup_i W_i=V$. For each $i$ let $\mu^\chi_i$
be the (finitely supported) uniform measure on ${\mathcal C}_q^\chi(W_i)$. Any
(weak) subsequential limit of the $\mu^\chi_i$'s (and by compactness there
must be at least one such) is a Gibbs measure. This fact was originally
proved, in a much more general context, by Dobrushin \cite{Dobrushin4}; see
e.g.\ \cite[Theorem 3.5]{BrightwellWinkler} for a simple proof in the present context.

A central concern in statistical physics
(again see \cite{Georgii2} for a thorough discussion)
is understanding when a particular system --- in our case the $q$-coloring model --- exhibits {\em phase coexistence}
(a.k.a. {\em phase transition}) on a given
infinite $G$, meaning that it admits more than one
Gibbs measure.
Actually, as we explain below, what we are really interested in is whether there are multiple Gibbs measures that are all substantial in an appropriate sense.

Our particular concern here is with
% the case $q=3$ and
$G={\mathbb Z}^d$, the usual nearest-neighbor graph on the $d$-dimensional integer lattice. This is a bipartite graph, with bipartition classes ${\mathcal E}$ (the {\em even} vertices, the set of lattice points the sum of whose coordinates is even) and ${\mathcal O}$ (the {\em odd} vertices). We will also use ${\mathcal E}$ and ${\mathcal O}$ for induced partition classes of subgraphs of ${\mathbb Z}^d$.
Intuition suggests that, {\em for large $d$}, the set of $3$-colorings of ${\mathbb Z}^d$ should
mainly consist of six classes, each identified by a
predominance of one of the colors on one of ${\mathcal E}$, ${\mathcal O}$, with the other two colors mainly assigned to the other partition class;  thus (again for large enough $d$) the set of Gibbs measures should include
six distinct measures
corresponding to these classes.
A well-known conjecture that this is the case goes back at least to
R. Koteck\'y circa 1985 (\cite{RK}; see e.g.\ \cite{Kotecky} for context),
although the explicit conjecture seems not to have appeared in print.

Our first main result verifies Koteck\'y's conjecture.
To state the result precisely, we set up some notation. Let $\chi(0,{\mathcal O}) \in {\mathcal C}_3$ be any $3$-coloring of ${\mathbb Z}^d$ satisfying $\chi|_{\mathcal O} \equiv 0$. For each $n \in {\mathbb N}$, let $W_n$ consist of the box $\{-n, \ldots, n\}^d$ together with all of the odd vertices of the box $\{-(n+1), \ldots, n+1\}^d$.
%Write $\mu_n$ for $\mu^{\chi(0,{\mathcal O})}_n$.
%%% DJG October 9 2012: commented out above line; removing notation ``\mu_n'' since it is so infrequently used
Let $v \in {\mathcal E}$ and $w \in {\mathcal O}$ be fixed vertices of ${\mathbb Z}^d$.
%\gbs{Should $v$, $w$ be in $W_n$, and if so does it matter where?}
Let $\mu^{(0,{\mathcal O})}$ be any subsequential limit of the $\mu^{\chi(0,{\mathcal O})}_n$'s.
%%% DJG Oct 9 2012: replaced a \mu_n above with \mu^{\chi(0,{\mathcal O})}_n
\begin{thm} \label{thm-influence1}
With notation as above,
$$
\mu^{(0,{\mathcal O})}(\sigma(v)=m)  \left\{
\begin{array}{lll}
\leq & e^{-\Omega(d)} & \mbox{if $m=0$} \\
\geq & 1/2 - e^{-\Omega(d)} & \mbox{if $m \in \{1,2\}$}
\end{array}
\right.
$$
and
$$
\mu^{(0,{\mathcal O})}(\sigma(w)=m)  \left\{
\begin{array}{lll}
\geq & 1/2 - e^{-\Omega(d)} & \mbox{if $m=0$} \\
\leq & 1/4 + e^{-\Omega(d)} & \mbox{if $m \in \{1,2\}$.}
\end{array}
\right.
$$
\end{thm}
This immediately implies that $\mu^{(0,{\mathcal O})}$ together with $\mu^{(1,{\mathcal O})}$, $\mu^{(2,{\mathcal O})}$, $\mu^{(0,{\mathcal E})}$, $\mu^{(1,{\mathcal E})}$ and $\mu^{(2,{\mathcal E})}$ (all defined in the obvious way) form a collection of six distinct Gibbs measures for all sufficiently large $d$.

It is possible for a Gibbs measure to be trivial. For example, if $\chi \in
{\mathcal C}_3({\mathbb Z}^2)$ is the mod 3 coloring (satisfying $\chi((x,y))
= x+y ~({\rm mod}~3)$) and $W_i$ is the $\ell_\infty$ ball of radius $i$, then
it is straightforward to check that the only coloring that agrees with $\chi$
off $W_i$ is $\chi$ itself, and so the $\mu^\chi_i$'s in this case have as
their unique limit the Gibbs measure with support $\{\chi\}$. (See e.g.\ \cite{BrightwellWinkler4} for other examples of such ``frozen'' Gibbs measures for the $q$-coloring model on the infinite regular tree.) These trivialities are avoided if we focus on Gibbs measures {\em of maximal entropy} (essentially measures with substantial support; see Section \ref{sec-entropyproof} for a precise definition). Koteck\'y's conjecture as originally told to us
\cite{RK} was that the 3-coloring model in high dimension
admits multiple Gibbs measures of maximal entropy.
\begin{thm} \label{thm-maxent}
The Gibbs measure $\mu^{(0,{\mathcal O})}$ constructed above is a measure of maximal entropy.
\end{thm}

At present our methods do not extend beyond $q=3$, but we strongly believe that the phenomenon of phase coexistence for the $q$-coloring model on ${\mathbb Z}^d$ occurs for all $q \geq 3$. A resolution of the following conjecture would be of great interest in both the statistical physics and discrete probability communities.
\begin{conj} \label{conj-q>3}
For all $q > 3$ and all sufficiently large $d=d(q)$, there is more than one Gibbs measure of maximal entropy for the $q$-coloring model on ${\mathbb Z}^d$.
\end{conj}
The natural expectation is that for odd $q$ there are at least $2{q \choose \lfloor q/2 \rfloor}$ such measures and for even $q$ at least ${q \choose q/2}$ such, corresponding to choices of a partition $[q]=A \cup B$ with $|A|=\lfloor q/2 \rfloor$ and a partition class of ${\mathbb Z}^d$ on which colors from $A$ are preferred. (Note that the issue here is only the analog of Theorem \ref{thm-influence1}; Theorem \ref{thm-maxent} extends without difficulty.) The analogous statement for proper $q$-colorings of the Hamming cube $\{0,1\}^d$ was proved in \cite{EngbersGalvin-tori}.

\medskip

In this paper we also consider the
problem of using Markov chains to sample uniformly at random from the set
${\mathcal C}_q(G)$, for finite $G$. Sampling and counting colorings of a graph are fundamental problems
in computer science and discrete mathematics. One approach is to design a Markov chain
whose stationary distribution is uniform over the set of
colorings of $G$. Then, starting from an arbitrary coloring and simulating a
random walk according to this chain for a sufficient number of
steps, we get a sample from a distribution which is close to uniform. The
number of steps required for the distribution to get close to uniform is
referred to as the {\em mixing time}   (see e.g.\ \cite{sinc}). The chain is called {\em rapidly mixing} if the mixing time is polynomial in $|V|$ (so
it converges quickly to stationarity); it is {\em torpidly mixing}
if its mixing time is super-polynomial in $|V|$
% DJG Oct 9 2012 --- changes $n$ to $|V|$
(so it converges
slowly). There has been a long history of studying mixing times of
various chains in the context of colorings (see e.g.\ \cite{ammb,
FriezeVigoda, gmp, hv, jer, lrs}).

A particular focus of this study has been on {\em Glauber dynamics}.
For $q$-colorings this is any single-site update Markov chain
that connects two colorings only if they differ on at most a single
vertex. The {\em Metropolis} chain ${\mathcal M}_q$ on state space ${\mathcal C}_q$
has transition probabilities $P_q(\chi_1,\chi_2)$, $\chi_1,\chi_2
\in {\mathcal C}_q,$ given by
$$
P_q(\chi_1,\chi_2) = \left\{
            \begin{array}{ll}
               0, & \mbox{ if $|\{v \in V:\chi_1(v)\neq \chi_2(v)\}| > 1$}; \\
               \frac{1}{q|V|},  &  \mbox{ if $|\{v \in V:\chi_1(v)\neq \chi_2(v)\}| = 1$}; \\
               1 - \sum_{\chi_1 \neq \chi_2' \in {\mathcal C}_q}  P_q(\chi_1,\chi_2') & \mbox{ if $\chi_1=\chi_2$.}
            \end{array}
         \right.
$$
% DJG Oct 9 2012 - changes a bunch of stray $k$'s above to $q$'s, and $q$ below to $[q]$
We may think of ${\mathcal M}_q$ dynamically as follows. From a $q$-coloring
$\chi$, choose a vertex $v$ uniformly from $V$ and a color $j$
uniformly from $[q]$.  Then recolor $v$ with color $j$
if the result is a (proper) $q$-coloring; otherwise stay at $\chi$.

When ${\mathcal M}_q$ is ergodic, its stationary distribution $\pi_q$ is
uniform over $q$-colorings. A series of recent papers has
shown that  ${\mathcal M}_q$ is rapidly mixing provided the number of colors
is sufficiently large compared to the maximum degree (see
\cite{FriezeVigoda} and the references therein). Substantially less
is known when the number of colors is small.  In fact, for $q$ small
it is NP-complete to decide whether a graph admits even one
$q$-coloring (see e.g.\ \cite{GareyJohnson}).

In this paper we consider the mixing rate of ${\mathcal M}_q$ on rectangular
regions of ${\mathbb Z}^d$. It is known \cite{lrs} that for $q \geq 3$, Glauber dynamics is ergodic (connects the
state space of $q$-colorings) on any such lattice region.
In ${\mathbb Z}^2$ much is known about the mixing rate of ${\mathcal M}_q$.
Randall and Tetali \cite{RandallTetali}, building on work of Luby et al.\ \cite{lrs}, showed that Glauber dynamics for sampling $3$-colorings
is rapidly mixing on any finite, simply-connected subregion of
${\mathbb Z}^2$ when the colors on the boundary of the region are fixed.
Goldberg et al.\ \cite{gmp} subsequently showed that the chain
remains fast on rectangular regions without this boundary
restriction. Substantially more is known when there are many colors:
Jerrum \cite{jer} showed that Glauber dynamics is rapidly mixing on
any graph satisfying $q \geq 2 \Delta$, where $q$ is the number of
colors and $\Delta$ is the maximum degree, thus showing Glauber
dynamics is fast on ${\mathbb Z}^2$ when $q \geq 8$.  It has since been shown
that it is fast for $q \geq 6$ \cite{ammb, bdg}.  Surprisingly, the
efficiency remains unresolved for $q=4$ or $5$.

% DJG October 9 2012 --- revised the following paragraph so that it more reflects reality
In higher dimensions much less is known when $q$ is small. The belief among
physicists working in the field is that Glauber dynamics on $3$-colorings is torpidly mixing
when the dimension $d$ of the cubic lattice is large enough (see e.g.\ the discussions in \cite{swk0,swk}), but there are no rigorous results. Here, we obtain the first such rigorous result by proving torpid mixing of the chain on cubic lattices with periodic boundary conditions.

%In higher dimensions much less is known when $q$ is small.
%Physicists have performed extensive numerical experiments (see e.g.\ \cite{fs,
%swk}) suggesting that Glauber dynamics on $3$-colorings is torpidly
%mixing
%%\gbs{for rectangular/square regions of $\mathbb Z^d$?}
%when the dimension
%%\gbs{$d$}
%of the cubic lattice is large enough. We
%prove this conjecture by studying the mixing
%time of the chain on cubic lattices with periodic boundary
%conditions.
%%\gbs{The result for the torus extends to rectangular/cubic regions
%%of $\mathbb Z^d$?}

Formally, we consider $3$-colorings of the even discrete
torus ${\mathbb Z}^d_n$. This is the graph on vertex set $[n]^d$ (with $n$ even) with edge set consisting of those pairs of
%%% DJG Oct 9 2012: changed $n$ to $[n]$ above
vertices that differ on exactly one coordinate and differ by $1$
(mod $n$) on that coordinate.
For a Markov chain ${\mathcal M}$ on ${\mathcal C}_3={\mathcal C}_3({\mathbb Z}^d_n)$ we denote by $\tau_{\mathcal M}$ the mixing time of
the chain (see Section \ref{sec-mixingproof} for a precise definition). We prove the following.
\begin{thm} \label{thm-glauber.slow.mixing}
There is a constant $d_0>0$ for which the following holds. For
$d\geq d_0$ and $n\geq 4$ even, the Glauber dynamics chain ${\mathcal M}_3$
on ${\mathcal C}_3$ satisfies
$$
\tau_{{\mathcal M}_3} \geq \exp \left\{\frac{n^{d-1}}{d^4\log^2 n}\right\}.
$$
\end{thm}
When $n=2$, ${\mathbb Z}^d_n$ becomes the Hamming cube $\{0,1\}^d$. Slow mixing of Glauber dynamics for sampling $3$-colorings was proved in this case in \cite{Galvin2}.
As with the case of phase coexistence, we strongly believe that torpid mixing holds for all $q>3$ as well, as long as the dimension is sufficiently high.
\begin{conj} \label{conj-q>3mixing}
For all $q > 3$, all even $n \geq 2$ and all sufficiently large $d=d(q)$, the mixing time of the Glauber dynamics chain ${\mathcal M}_q$
on ${\mathcal C}_q$ is (essentially) exponential in $n^{d-1}$.
\end{conj}

Our techniques actually apply to a more general class of chains. A
Markov chain ${\mathcal M}$ on state space ${\mathcal C}_3$ is said to be {\em $\rho$-local} if, in each step of the chain, at most $\rho|V|$ vertices have
their colors changed; that is, if
$$
P_{\mathcal M}(\chi_1,\chi_2) \neq 0 ~ \Rightarrow ~ |\{v \in V : \chi_1(v)\neq
\chi_2(v)\}| \leq \rho |V|.
$$
These types of chains were introduced in \cite{DyerFriezeJerrum},
where the terminology {\em $\rho|V|$-cautious} was employed. We
prove the following, which easily implies
Theorem~\ref{thm-glauber.slow.mixing}.
\begin{thm} \label{thm-rho.local.slow.mixing}
Fix $\rho>0$ satisfying $H(\rho)+\rho < 1$. There is a constant
$d_0=d_0(\rho)>0$ for which the following holds. For $d\geq d_0$ and
$n\geq 4$ even, if ${\mathcal M}$ is an ergodic $\rho$-local Markov chain on
${\mathcal C}_3$ with uniform stationary distribution then
$$
\tau_{{\mathcal M}} \geq \exp \left\{\frac{n^{d-1}}{d^4\log^2 n}\right\}.
$$
\end{thm}
Here $H(x)=-x\log x -(1-x)\log (1-x)$ is the usual binary entropy
function. Note that all $\rho \leq 0.22$ satisfy $H(\rho)+\rho < 1$.

\medskip

We show phase transition using a Peierls argument, to be discussed in detail in Section \ref{sec-preview}.
We show torpid mixing via a
conductance argument by identifying a cut in the state space
requiring exponential time to cross. For both results, our work
builds heavily on technical machinery introduced by Galvin and Kahn
\cite{GalvinKahn} showing that the hard-core (independent set) model on ${\mathbb Z}^d$ exhibits phase transition for some values of the density parameter $\lambda$ that go to zero as the dimension grows. Specifically, for $\lambda>0$, choose
${\mathbb I}$ from ${\cal I}(\Lambda_n)$ (the set of independent
sets of the box $\Lambda_n=\{-n, \ldots, n\}$) with $\Pr({\mathbb I}=I) \propto \lambda^{|I|}$.
Galvin and Kahn showed that for $\lambda > Cd^{-1/4}\log^{3/4}d$
(for some constant $C$) and fixed $v \in {\mathcal E}$,
$$
\lim_{n\rightarrow\infty}{\mathbb P}(v\in{\mathbb I} ~|~{\mathbb I}  \supseteq
\partial_{\rm int} \Lambda_n  \cap {\cal E})
 \geq \frac{(1+o(1))\lambda}{1+\lambda}
$$
whereas
$$
\lim_{n\rightarrow\infty}{\mathbb P}\left(v\in{\mathbb I}~|~{\mathbb
I}\supseteq
\partial_{\rm int} \Lambda_n\cap {\cal O}\right) \leq (1+\lambda)^{-2d(1-o(1))}
$$
where $\partial_{\rm int} \Lambda_n$ is the set of vertices in $\Lambda_n$ that are adjacent (in ${\mathbb Z}^d$) to something outside $\Lambda_n$.
In other words, the influence of the boundary on the center of a
large box persists as the boundary recedes.

Neither the results of \cite{GalvinKahn} (showing phase coexistence for the hard-core model on ${\mathbb Z}^d$) nor
Theorem \ref{thm-influence1} (concerning $3$-colorings of ${\mathbb Z}^d$) directly imply anything about the behavior of Markov chains on
finite lattice regions. However, they do suggest that in the finite
setting, typical configurations fall into the distinct classes
described in stationarity and that local Markov chains will be unlikely to move
between these classes; the remaining configurations are expected to
have negligible weight for large lattice regions, even when they are
finite.

Galvin \cite{Galvin} extended the results of \cite{GalvinKahn},
showing that in sufficiently high dimension, Glauber dynamics on
independent sets mixes torpidly in rectangular regions of ${\mathbb
Z}^d$ with periodic boundary conditions. Similar results were known
previously about independent sets; however, one significant new
contribution of \cite{Galvin} was showing that as $d$ increases, the
critical $\lambda$ above which Glauber dynamics mixes torpidly tends
to $0$. In particular, there is some dimension $d_0$ such that for
all $d \geq d_0$, Glauber dynamics will be torpid on ${\mathbb Z}^d$ when
$\lambda=1$. This turns out to be the crucial new ingredient
allowing us to rigorously verify phase transition and torpid mixing for
$3$-colorings in high dimensions, as there turns out to be a close
connection between the independent set model at $\lambda=1$ and the
$3$-coloring model. % Note that u
Unlike many statistical physics
models, the $3$-coloring model does not come equipped with a parameter such as $\lambda$
that can be tweaked to establish desired bounds; this makes the
proofs here significantly more delicate than the usual phase-transition and torpid-mixing
arguments.

\medskip

The rest of the paper is laid out as follows. In Section \ref{sec-influenceproof} we give the proof of Theorem \ref{thm-influence1} (phase transition), modulo one of our two main technical lemmas, Lemma \ref{lem-volume.bounds.from.gk-pt}. This section also provides an overview of our proof strategy (Section \ref{sec-preview}). In Section \ref{sec-mixingproof}, we give the proof of Theorem \ref{thm-rho.local.slow.mixing} (torpid mixing), modulo the second main technical lemma, Lemma \ref{lem-volume.bounds.from.gk}. Section \ref{sec-lemmaproofs} provides the proofs of Lemmas \ref{lem-volume.bounds.from.gk-pt} and \ref{lem-volume.bounds.from.gk}, while in Section \ref{sec-entropyproof} we prove Theorem \ref{thm-maxent} (measures of maximal entropy).

\medskip

The original aim of this work was to prove Theorem \ref{thm-influence1} (phase coexistence).
We achieved this at a 2002 Newton Institute programme,\footnote{Isaac Newton Institute for Mathematical Sciences programme on
Computation, Combinatorics and Probability, 29 Jul-- 20 Dec 2002,
\texttt{http://www.newton.ac.uk/programmes/CMP/}.}
and the second author discussed the result in talks
and in communications with R. Koteck\'y and others.
We then noticed that with some additional work we could obtain a proof of Theorem \ref{thm-rho.local.slow.mixing} (torpid mixing). This result was presented by the first and third authors in \cite{GalvinRandall},
which also includes the first mention of Theorem \ref{thm-influence1} in print.
During preparation
of the present manuscript we learned from Ron Peled
of his recent \cite{Peled},
whose main result contains Theorem \ref{thm-influence1} (of which he heard from Koteck\'y only after proving
his result \cite{RP}).
Though similar in spirit, the approach of \cite{Peled}, which exploits a correspondence
between  colorings and height functions, is different from the present argument, which
stays within the world of colorings.

\section{Proof of Theorem \ref{thm-influence1}} \label{sec-influenceproof}

In this section we show that the $3$-coloring model on ${\mathbb Z}^d$ admits multiple Gibbs measures for all sufficiently large $d$ (Theorem \ref{thm-influence1}).

\subsection{Some notation}

Let $\Sigma=(V,E)$ be a bipartite graph with bipartition classes ${\mathcal E}$ and ${\mathcal O}$.
For $X \subseteq V$, write $\nabla(X)$ for the set of edges in $E$
that have one end in $X$ and one end outside $X$; $\overline{X}$ for
$V \setminus X$; $\partial_{\rm int}X$ for the set of vertices in $X$
that are adjacent to something outside $X$; $\partial_{\rm ext}X$ for
the set of vertices outside $X$ that are adjacent to something in
$X$; $X^+$ for $X \cup
\partial_{\rm ext} X$; $X^{\mathcal E}$ for $X \cap {\mathcal E}$ and $X^{\mathcal O}$ for $X \cap
{\mathcal O}$. Further, for $x \in V$ set $\partial x=\partial_{\rm ext}\{x\}$.
We abuse notation slightly, identifying
sets of vertices of $V$ and the subgraphs they induce.

\subsection{Finitizing Theorem \ref{thm-influence1}}

Theorem \ref{thm-influence1} may be finitized as follows. Set $\Lambda=\Lambda_n=\{-n, \ldots, n\}^d$. (Throughout, $n$ will be fixed, so we drop the dependence in the notation.)
Set
$$
{\mathcal C}_3^{\mathcal O} = \{\chi \in {\mathcal C}_3(\Lambda) : \chi|_{\partial_{\rm int} \Lambda \cap {\mathcal O}}
\equiv 0\}
$$
and for $v_0 \in (\Lambda \setminus \partial_{\rm int} \Lambda) \cap {\mathcal E}$ set
$$
{\mathcal C}_3^{\mathcal O}(v_0) = \{\chi \in {\mathcal C}_3^{\mathcal O} : \chi(v_0)=0\}.
$$
In other words, ${\mathcal C}_3^{\mathcal O}$ is the subset of (proper) $3$-colorings of $\Lambda$ in which
all of the odd vertices on the boundary get color $0$, while
${\mathcal C}_3^{\mathcal O}(v_0)$ is the set of those colorings in which even $v_0$ also
gets color $0$. We prove the following.
\begin{thm} \label{thm-influence2}
For all $n$,
\begin{equation} \label{inq-odd.boundary.on.even}
\frac{|{\mathcal C}_3^{\mathcal O}(v_0)|}{|{\mathcal C}_3^{\mathcal O}|} \leq e^{-\Omega(d)}
\end{equation}
as $d \rightarrow \infty$ (with the implicit constant independent of $n$).
\end{thm}
With some extra work we could replace $e^{-\Omega(d)}$ here with $2^{-2d(1-o(1))}$. This would require dealing more carefully with small $c_0$ in Lemma \ref{lem-volume.bounds.from.gk-pt}, and to simplify the presentation we chose not to do this. The interested reader may consult \cite{GalvinKahn} (and in particular the end of Section 2.13 of that reference) for the approach.

Theorem \ref{thm-influence2} implies Theorem \ref{thm-influence1}. Indeed, let $\mu^{(0,{\mathcal O})}$ be any subsequential limit of the $\mu^{\chi(0,{\mathcal O})}_n$'s, where the notation
%%% DJG Oct 9 2012: replaced a \mu_n above with \mu^{\chi(0,{\mathcal O})}_n
is as in the discussion before the statement of Theorem \ref{thm-influence1}. From Theorem \ref{thm-influence2} we immediately have $\mu^{(0,{\mathcal O})}(\chi(v)=0) \leq e^{-\Omega(d)}$ and so, by symmetry, $\mu^{(0,{\mathcal O})}(\chi(v)=1)=\mu^{(0,{\mathcal O})}(\chi(v)=2) \geq 1/2 - e^{-\Omega(d)}$. A second application of Theorem \ref{thm-influence2} (together with a union bound) shows that, for $Q=\{\chi(v') \neq 0~\forall v' \sim w\}$, $\mu^{(0,{\mathcal O})}(Q) = 1 - e^{-\Omega(d)}$, whence
$$
\mu^{(0,{\mathcal O})}(\chi(w)=0) = \mu^{(0,{\mathcal O})}(Q)\mu^{(0,{\mathcal O})}(\chi(w)=0|Q) \geq (1-e^{-\Omega(d)})/2;
$$
and then symmetry gives the second inequality in Theorem \ref{thm-influence1}.

\subsection{Preview} \label{sec-preview}

For a generic $\chi \in {\mathcal C}_3^{\mathcal O}(v_0)$ there is a region of
$\Lambda$ around $v_0$ consisting predominantly of even vertices colored $0$
together with their neighbors, and a region around $\partial_{\rm int} \Lambda$ consisting of odd
vertices colored $0$ together with their neighbors. These regions
are separated by a
two-layer $0$-free moat or {\em cutset}. In Section \ref{sec-cutsets} we describe a procedure that associates a
particular such cutset with each $\chi  \in {\mathcal C}_3^{\mathcal O}(v_0)$. Our main technical result,
Lemma \ref{lem-volume.bounds.from.gk-pt}, asserts that for each possible cutset size $c$, the probability that the cutset associated with a uniformly chosen coloring has size $c$ is
exponentially small in $c$. This lemma is
presented in Section \ref{sec-main-tech-lemma}, where it is also used to derive Theorem \ref{thm-influence2}.

We use a variant of the Peierls argument (originally presented in \cite{Peierls}) to prove Lemma
\ref{lem-volume.bounds.from.gk-pt}. By carefully modifying $\chi
\in {\mathcal C}_3^{\mathcal O}(v_0)$ inside its cutset,
we can exploit the fact that the cutset is $0$-free to map $\chi$
to a set $\varphi(\chi)$ of many different $\chi' \in {\mathcal C}_3^{\mathcal O}$.
If the $\varphi(\chi)$'s were disjoint for distinct $\chi$'s, we
would be done, having shown that there are many more
$3$-colorings in ${\mathcal C}_3^{\mathcal O}$ than in ${\mathcal C}_3^{\mathcal O}(v_0)$.
To control the possible overlap,
we define a flow $\nu:{\mathcal C}_3^{\mathcal O}(v_0) \times{\mathcal C}_3^{\mathcal O}\rightarrow
[0,\infty)$ supported on pairs $(\chi,\chi')$ with $\chi' \in
\varphi(\chi)$ in such a way that the flow out of each $\chi \in
{\mathcal C}_3^{\mathcal O}(v_0)$ is $1$. Any uniform bound we can obtain on the flow
into elements of ${\mathcal C}_3^{\mathcal O}$
is then easily seen to be a bound on $|{\mathcal C}_3^{\mathcal O}(v_0)|/|{\mathcal C}_3^{\mathcal O}|$. We
define the flow via a notion of approximation modified from
\cite{GalvinKahn}. To each cutset $\gamma$ we associate a set
$A(\gamma)$ that approximates the interior of $\gamma$ in a precise
sense, in such a way that as we run over all possible $\gamma$, the
total number of approximate sets used is small. Then for each $\chi'
\in {\mathcal C}_3^{\mathcal O}$ and each approximation $A$, we consider the set of those $\chi \in {\mathcal C}_3^{\mathcal O}(v_0)$ with
$\chi' \in \varphi(\chi)$ and with $A$ the approximation to
$\gamma$. We define the flow so that if this set is large, then
$\nu(\chi,\chi')$ is small for each $\chi$ in the set. In this way
we control the flow into $\chi'$ corresponding to each
approximation $A$; since the total number of
approximations is small, we control the total flow into $\chi'$. In
the language of statistical physics,
this approximation scheme
is a {\em course-graining} argument. The details appear in
Section \ref{sec-lemmaproofs}.

The main result of \cite{GalvinKahn} is proved
along similar lines to those described above. One of the
difficulties we encounter in moving from these arguments on
independent sets to arguments on colorings is that of finding an
analogous way of modifying a coloring inside a cutset in order to
exploit the fact that it is $0$-free. The beginning of Section
\ref{sec-lemmaproofs} (in particular Claims \ref{claim-shift.works} and
\ref{claim-unique.reconstruction}) describes an appropriate
modification that has all the properties we desire.

\subsection{Cutsets} \label{sec-cutsets}

We now describe a way of associating with each $\chi \in {\mathcal C}_3^{\mathcal O}(v_0)$
a minimal edge cutset, following an approach of
\cite{BorgsChayesFriezeKimTetaliVigodaVu} and \cite{Galvin}. (An alternate construction is given in \cite{GalvinKahn}. The present construction is
perhaps more transparent.)

Given $\chi \in {\mathcal C}_3^{\mathcal O}(v_0)$ set $I=I(\chi)=\chi^{-1}(0)$.
Note that $I$ is an independent set (a set of vertices no two
of which are adjacent). Let $R$ be the component of $(I^{\mathcal E})^+$ that includes $v_0$. Let $C$ be the component of $\overline{R}$, that includes $\partial_{\rm int} \Lambda$. Set $\gamma =
\gamma(\chi)=\nabla(C)$ and $W=W(\chi)=\overline{C}$. Evidently
$C$ is connected, and $W$ consists of $R$, which is connected,
together with a number of other components of $\overline{R}$, each
of which are joined to $R$; so $W$ is also connected.
It follows that $\gamma$ is a minimal edge-cutset in $\Lambda$, separating $v_0$ from $\partial_{\rm int} \Lambda$. Note that $\gamma$ depends
only on the independent set $I$. Note also that the vertex set of $\gamma$ is $\partial_{\rm int} W \cup \partial_{\rm ext} W$. We write $|\gamma|$ for the size (number of edges) of $\gamma$.

The next lemma summarizes the
properties of $\gamma$ that we will draw upon in what follows; having
established these properties we will not subsequently refer to the
details of the construction. For the most part these
properties will not be used directly, but will be referred to
to validate the applications of various results from
\cite{GalvinKahn}.

\begin{lemma} \label{lem-properties.of.contours}
For each $\chi \in {\mathcal C}_3^{\mathcal O}(v_0)$ we have the following.
\begin{equation} \label{contour.prop.1}
v_0 \in W ~~\mbox{and}~~ \partial_{\rm int} \Lambda \cap W = \emptyset;
\end{equation}
\begin{equation} \label{contour.prop.2}
\partial_{\rm int}W \subseteq {\mathcal O}~~\mbox{and}~~\partial_{\rm ext}W \subseteq
{\mathcal E};
\end{equation}
\begin{equation} \label{contour.prop.3}
\partial_{\rm int}W \cap I =
\emptyset~~\mbox{and}~~\partial_{\rm ext}W \cap I = \emptyset;
\end{equation}
\begin{equation} \label{contour.prop.4}  %QQ (label used?)
\forall v \in \partial_{\rm int}W, ~\partial  v \cap W \cap I \neq
\emptyset;
\end{equation}
\begin{equation} \label{contour.prop.5}
W^{\mathcal O} = \partial_{\rm ext} W^{\mathcal E}~~\mbox{and}~~W^{\mathcal E}=\left\{y \in
{\mathcal E}:\partial y \subseteq W^{\mathcal O}\right\}
\end{equation}
and
\begin{equation} \label{contour.prop.8}
\mbox{for large enough $d$, ~~$|\gamma| \geq \max\{|W|^{1-1/d}, d^2\}$}.
\end{equation}
\end{lemma}

\smallskip

\noindent {\em Proof}: That $v_0 \in W$ and $\partial_{\rm int} \Lambda \cap W =
\emptyset$ is clear.

Properties (\ref{contour.prop.2}), (\ref{contour.prop.3}), (\ref{contour.prop.4}) and (\ref{contour.prop.5}) are also easily verified; see
\cite[Lemma 3.3]{Galvin} or \cite[Proposition 2.6]{GalvinKahn} for detailed proofs.

The isoperimetric inequality of Bollob\'as and Leader \cite[Theorem 3]{BollobasLeader2} says that if $W \subseteq \Lambda$ satisfies $|W| \leq n^d/2$ then $|\nabla(W)| \geq |W|^{1-1/d}$. Since $W \cap \partial_{\rm int} \Lambda = \emptyset$ we may apply this (perhaps with $W$ viewed as a subset of a larger $\Lambda$) to conclude that $|\gamma| \geq |W|^{1-1/d}$. For the second inequality in (\ref{contour.prop.8}), note that by (\ref{contour.prop.5}) we have $|\gamma|=2d(|W^{\mathcal O}|-|W^{\mathcal E}|)$. In \cite[Lemma 2.13]{GalvinKahn} it is shown that if $A, B \subseteq \Lambda$ satisfy $A \subseteq {\mathcal E}$, $B
\subseteq {\mathcal O}$, $B = \partial_{\rm ext} A$, $A=\{v \in {\mathcal E}: \partial v
\subseteq B\}$, $(A \cup B) \cap \partial \Lambda = \emptyset$, and $|B| < d^{O(1)}$ then $|B|-|A| \geq |B|(1-O(1/d))$.
By (\ref{contour.prop.1}) and (\ref{contour.prop.5}), $W^{\mathcal O}$ and $W^{\mathcal E}$ satisfy these conditions, and so noting that $|W^{\mathcal O}|\geq 2d$, we get $|\gamma| \geq 2d^2(1-o(1)) \geq d^2$ for $|W| \leq d^{O(1)}$; the first inequality in (\ref{contour.prop.8}) implies the second for all larger $|W|$.
\qed

\medskip

The cutsets also satisfy a connectivity property (specifically, that $\partial_{\rm int} W \cup \partial_{\rm ext} W$ induces a connected graph). We will not use this property explicitly in the sequel; it is an important ingredient in the proof of Lemma \ref{lem-k.comp.approx} (a combination of results from \cite{Galvin} and \cite{GalvinKahn}) where it serves to bound the number of cutsets of a given size that use a given edge.

\subsection{The main lemma for phase transition} \label{sec-main-tech-lemma}

For $c_0 \in {\mathbb N}$ set
$$
{\mathcal W}(c_0,v_0) = \left\{\gamma : |\gamma|=c_0,~\gamma = \gamma(\chi) ~\mbox{for some}~ \chi \in {\mathcal C}_3^{\mathcal O}(v_0)\right\}
$$
and set ${\mathcal W} = \cup_{c_0} {\mathcal W}(c_0,v_0)$. Set
${\mathcal C}_3^{\mathcal O}(c_0,v_0) = \left\{\chi \in {\mathcal C}_3^{\mathcal O}(v_0) : |\gamma(\chi)|=c_0\right\}.$
The main technical lemma we need to prove phase transition is the following.
\begin{lemma} \label{lem-volume.bounds.from.gk-pt}
There are constants $C, d_0>0$ such that the following holds. For
all $d \geq d_0$, $n$ and $c_0$,
$$
\frac{|{\mathcal C}_3^{\mathcal O}(c_0,v_0)|}{|{\mathcal C}_3^{\mathcal O}|} \leq \exp\left\{-\frac{Cc_0}{d}\right\}.
$$
\end{lemma}
We give the proof in Section \ref{sec-lemmaproofs}.

From Lemma \ref{lem-volume.bounds.from.gk-pt}, we easily obtain Theorem \ref{thm-influence2}. Indeed, for all $n$ and $d\geq d_0$ we have (using (\ref{contour.prop.8}) for the restriction on $c_0$)
\begin{eqnarray*}
|{\mathcal C}_3^{\mathcal O}(v_0)| & \leq & \sum_{c_0 \geq d^2} |{\mathcal C}_3^{\mathcal O}(c_0,v_0)| \\
& \leq & \sum_{c_0 \geq d^2} \exp\left\{-\frac{Cc_0}{d}\right\}|{\mathcal C}_3^{\mathcal O}| \\
& \leq & e^{-\Omega(d)} |{\mathcal C}_3^{\mathcal O}|.
\end{eqnarray*}

\section{Proof of Theorem \ref{thm-rho.local.slow.mixing}} \label{sec-mixingproof}

The aim of this section is to show that in the finite setting of the discrete torus, Glauber dynamics for sampling from $3$-colorings mixes torpidly (Theorem \ref{thm-rho.local.slow.mixing}). We begin by formalizing some definitions. Given an ergodic Markov chain
${\cal M}$ on state space $\Omega$ with stationary
distribution $\pi$, let $P^t(x,
\cdot)$ be the distribution of the chain at time $t$ given that it
started in state $x$.
The {\em mixing time} $\tau_{{\cal M}}$ of ${\mathcal M}$ is defined to be
$$
\tau_{{\mathcal M}}=\min \left\{t_0: \max_{x \in \Omega} \frac{1}{2} \sum_{y \in \Omega}
|P^t(x,y) - \pi(y)|  \leq
\frac{1}{e} ~~~ \forall t>t_0\right\}.
$$

We prove Theorem \ref{thm-rho.local.slow.mixing} via a well-known
{\em conductance argument} \cite{JerrumSinclair, LawlerSokal, Thomas}, using
a form of the argument derived in \cite{DyerFriezeJerrum}. Let
$A \subseteq \Omega$ and $M \subseteq \Omega \setminus A$ satisfy $\pi(A)
\leq 1/2$ and $\omega_1 \in A, \omega_2 \in \Omega \setminus (A \cup M)
\Rightarrow P(\omega_1, \omega_2) =0$. Then from \cite{DyerFriezeJerrum}
we have
\begin{equation} \label{conductance_bound}
\tau_{\mathcal M} \geq \frac{\pi(A)}{8\pi(M)}.
\end{equation}

Let us return to the setup of Theorem
\ref{thm-rho.local.slow.mixing}. For even $n$, ${\mathbb Z}^d_n$ is
bipartite with partition classes ${\mathcal E}$ (consisting of those vertices
the sum of whose coordinates is even) and ${\mathcal O}$. We will show that most
3-colorings have an imbalance whereby the vertices colored $0$ lie
either predominantly in ${\mathcal E}$ or predominantly in ${\mathcal O}$, and those that are roughly balanced are highly unlikely in stationarity.
Accordingly let us define the set of balanced $3$-colorings by
$$
{\mathcal C}_3^{b,\rho} = \{\chi \in {\mathcal C}_3:\left||\chi^{-1}(0)\cap
{\mathcal E}|\!-\!|\chi^{-1}(0)\cap {\mathcal O}|\right| \leq \rho n^d/2\}
$$
and let
$$
{\mathcal C}_3^{{\mathcal E},\rho} = \{\chi \in {\mathcal C}_3:|\chi^{-1}(0)\cap
{\mathcal E}|-|\chi^{-1}(0)\cap {\mathcal O}| > \rho n^d/2\}.
$$
By symmetry,
$\pi_3({\mathcal C}_3^{{\mathcal E},\rho}) \leq 1/2$ (recall that $\pi_3$ is uniform distribution). Notice that since ${\mathcal M}$
updates at most $\rho n^d$ vertices in each step, we have that
if $\chi_1 \in {\mathcal C}_3^{{\mathcal E},\rho}$ and $\chi_2 \in {\mathcal C}_3 \setminus
({\mathcal C}_3^{{\mathcal E},\rho} \cup {\mathcal C}_3^{b,\rho})$ then $P_{\mathcal M}(\chi_1,\chi_2) =
0$.
Therefore, by (\ref{conductance_bound}),
$$
\tau_{\mathcal M} \geq \frac{\pi_3({\mathcal C}_3^{{\mathcal E},\rho})}{8\pi_3({\mathcal C}_3^{b,\rho})}
\geq \frac{1-\pi_3({\mathcal C}_3^{{\mathcal E},\rho})}{16\pi_3({\mathcal C}_3^{b,\rho})},
$$
and so Theorem \ref{thm-rho.local.slow.mixing} follows from the
following critical theorem.
\begin{thm} \label{thm-main-mixing}
Fix $\rho>0$ satisfying $H(\rho)+\rho < 1$. There is a constant
$d_0=d_0(\rho)>0$ for which the following holds. For $d \geq d_0$
and $n \geq 4$ even,
$$
\pi_3({\mathcal C}_3^{b,\rho}) \leq \exp\left\{\frac{-2n^{d-1}}{d^4\log^2
n}\right\}.
$$
\end{thm}

\subsection{Cutsets revisited} \label{sec-cutsets2}

One difficulty we have to overcome
in moving from a Gibbs measure argument to a torpid mixing argument
is that of going from bounding the probability of a configuration
having a single cutset to bounding the probability of it having an
ensemble of cutsets. Another difficulty is that the cutsets we
consider in these ensembles can be topologically more complex than
the connected cutsets that are considered in the phase transition result. In part,
both of these difficulties are dealt with by the machinery developed
in \cite{Galvin}.

We begin by describing a way of associating with each $\chi \in {\mathcal C}_3^{b,\rho}$
a collection of minimal edge cutsets, extending the process described in Section \ref{sec-cutsets}.

For $\chi \in {\mathcal C}_3^{b,\rho}$ set $I=I(\chi)=\chi^{-1}(0)$. Given a component $R$ of $(I^{\mathcal E})^+$ or
$(I^{\mathcal O})^+$ and a component $C$ of $\overline{R}$, set $\gamma =
\gamma(R,C,\chi)=\nabla(C)$ and $W=W(R,C,\chi)=\overline{C}$. As in Section \ref{sec-cutsets}, $\gamma$ is a minimal edge-cutset in ${\mathbb Z}^d_n$.
Define ${\rm int} \gamma$, the {\em interior} of
$\gamma$, to be the smaller of $C,W$ (if $|W|=|C|$, take ${\rm int}
\gamma = W$).

The collection of cutsets associated to $\chi$ depend only on the independent set $I$, and coincide exactly with the
cutsets associated to an independent set in \cite{Galvin}. We may
therefore apply the machinery developed in \cite{Galvin} for
independent set cutsets in the present setting. In particular, from
\cite[Lemmas 3.1 and 3.2]{Galvin} we know that for each $\chi \in
{\mathcal C}_3$ there is a subset $\Gamma(\chi)$ of the collection of cutsets associated to $\chi$ that either satisfies
\begin{equation} \label{cutset.conditions}
\begin{array}{c}
\mbox{for all $\gamma \in \Gamma(\chi)$, ${\rm int} \gamma = W$, for all $\gamma,\gamma' \in \Gamma(\chi)$ with $\gamma \neq \gamma'$, ${\rm int} \gamma \cap {\rm int} \gamma' = \emptyset$,} \\
\mbox{and for all $\gamma \in \Gamma(\chi)$, $R$ is a component of $(I^{\mathcal E})^+$ and $I^{\mathcal E} \subseteq \cup_{\gamma \in \Gamma(\chi)} {\rm int} \gamma$},
\end{array}
\end{equation}
or the analogue of (\ref{cutset.conditions}) with ${\mathcal E}$ replaced by
${\mathcal O}$. Set ${{\mathcal C}}_3^{\rm even} = \{\chi \in {\mathcal C}_3:\chi~\mbox{satisfies
(\ref{cutset.conditions})}\}$. From here on whenever $\chi \in
{{\mathcal C}}_3^{\rm even}$ is given we assume that $I$ is its associated
independent set and that $\Gamma(\chi)$ is a particular collection of
cutsets associated with $\chi$ and satisfying
(\ref{cutset.conditions}).

The cutsets that we have constructed here have many properties in common with those constructed in Section \ref{sec-cutsets}; in particular,
each $\gamma \in \Gamma(\chi)$ satisfies (\ref{contour.prop.2}), (\ref{contour.prop.3}), (\ref{contour.prop.4}), (\ref{contour.prop.5}) and (\ref{contour.prop.8}). The proof of (\ref{contour.prop.8}) appeals to \cite[Theorem 8]{BollobasLeader2} instead of \cite[Theorem 3]{BollobasLeader2} and uses the fact that for large enough $n$ and for $|B|=d^{O(1)}$ we may apply \cite[Lemma 2.13]{GalvinKahn} in the setting of the torus without modification.

The cutsets in $\Gamma(\chi)$ also satisfy a connectivity property, although because the torus is topologically more complex than ${\mathbb Z}^d$ the connectivity property is more involved. In \cite[Lemma 3.4]{Galvin} it is shown that each $\gamma \in \Gamma(\chi)$ is either connected in the dual of the torus (the graph on the edges of the torus in which two edges are adjacent if there is a $4$-cycle including both of them) or has at least $n^{d-1}$ edges in each component. As in the case of phase transition, this property is important in the proof of Lemma \ref{lem-k.comp.approx}, but since we take this lemma directly from \cite{Galvin} we do not give further details here.

\subsection{The main lemma for torpid mixing} \label{subsec-results.from.GalvinKahn}

For $c \in {\mathbb N}$ and $v \in V$ set
$$
{\mathcal W}(c,v) = \left\{\gamma : |\gamma|=c,~ \gamma \in \Gamma(\chi)~\mbox{for some}~\chi \in {\mathcal C}_3^{\rm even},~\mbox{and}~ v \in ({\rm int} \gamma)^{\mathcal E}
\right\}
$$
and set ${\mathcal W} = \cup_{c,v} {\mathcal W}(c,v)$. A {\em profile} of a collection
$\{\gamma_0, \ldots, \gamma_\ell\} \subseteq {\mathcal W}$ is a vector
$\ul{p}=(c_0 ,v_0, \ldots, c_\ell, v_\ell)$ with $\gamma_i \in
{\mathcal W}(c_i,v_i)$ for all $i$. Given a profile $\ul{p}$ set
$$
{\mathcal C}_3(\ul{p}) = \left\{\chi \in {\mathcal C}_3^{\rm even} : \Gamma(\chi)~\mbox{contains a subset with profile}~ \ul{p} \right\}.
$$
Our main lemma (c.f. \cite[Lemma 3.5]{Galvin})
is the following.
\begin{lemma} \label{lem-volume.bounds.from.gk}
There are constants $C, d_0>0$ such that the following holds. For
all even $n\geq 4$ and $d \geq d_0$, and all profiles $\ul{p}$ as above,
$$
\pi_3({\mathcal C}_3(\ul{p})) \leq \exp\left\{-\frac{C\sum_{i=0}^\ell c_i}{d}
\right\}.
$$
\end{lemma}

\subsection{Proof of Theorem \ref{thm-main-mixing}} \label{subsec-proof_of_thm}

We will prove Lemma \ref{lem-volume.bounds.from.gk} in Section \ref{sec-lemmaproofs}. Here, we derive Theorem \ref{thm-main-mixing} from it. Throughout we assume that the conditions of Theorem \ref{thm-main-mixing} and Lemma
\ref{lem-volume.bounds.from.gk} are satisfied (with $d_0$
sufficiently large to support our assertions).

We begin with an easy count that dispenses with colorings where
$|I(\chi)|$ is small. Set
$$
{\mathcal C}_3^{\rm small} = \left\{\chi \in
{\mathcal C}_3^{b,\rho}:\min\{|I^{\mathcal E}|,|I^{\mathcal O}|\} \leq \frac{n^d}{4d^{1/2}}\right\}.
$$
\begin{lemma} \label{lem-bounding_small} $\pi_3({\mathcal C}_3^{\rm small}) \leq
\exp\left\{-\Omega(n^d)\right\}$.
\end{lemma}

\smallskip

\noindent {\em Proof}:
For any $A \subseteq {\mathcal E}$ and $B \subseteq {\mathcal O}$, let
${\rm comp}(A,B)$ be the number of components in
$V \setminus (A \cup B \cup \partial^\star A \cup \partial^\star B)$, where
for $T \subseteq {\mathcal E}$ (or ${\mathcal O}$),
$$
\partial^\star T = \{x \in \partial_{\rm ext} T : \partial x \subseteq T\}~(=\{x \in V :
\partial x \subseteq T\}).
$$
We begin by noting that by ${\mathcal E}$-${\mathcal O}$
symmetry
\begin{equation} \label{inq-small.count}
|{\mathcal C}_3^{\rm small}| \leq 2 \sum \exp_2\left\{|\partial^\star A|+|\partial^\star
B|+{\rm comp}(A,B)\right\},
\end{equation}
where the sum is over all pairs $A \subseteq {\mathcal E}$, $B \subseteq {\mathcal O}$
with no edges between $A$ and $B$ and satisfying $|A|\leq
n^d/4d^{1/2}$ and $|B| \leq (\rho + 1/2d^{1/2})n^d/2$.
Indeed, once
we have specified that the set of vertices colored $0$ is $A \cup
B$, we have a free choice between $1$ and $2$ for the color at $x
\in  \partial^\star A \cup
\partial^\star B$, and we also have a free choice between the two
possible colorings of each component of $V \setminus (A \cup B \cup
\partial^\star A \cup \partial^\star B)$.

A key observation is the following.  For $A$ and $B$ contributing to the sum in
(\ref{inq-small.count}),
\begin{equation} \label{inq-number.of.remaining.comps}
{\rm comp}(A,B) \leq \frac{n^d}{2d}.
\end{equation}
To see this, let $C$ be a component of $V \setminus (A \cup B)$. If
$C=\{v\}$ consists of a single vertex, then (depending on the parity
of $v$) we have either $\partial v \subseteq A$ or $\partial v
\subseteq B$ and so $v \in \partial^\star A \cup \partial^\star B$.
Otherwise, let $vw$ be an edge of $C$ with $v \in {\mathcal E}$ (and so $w
\in {\mathcal O}$). If $v$ has $k$ edges to $B$ and $u$ has $\ell$ to $A$,
then (since there are no edges from $A$ to $B$)
we have $(k-1)+(\ell-1)\leq 2d-2$ or $k+\ell \leq 2d$. (Here we are
using that in ${\mathbb Z}^d_n$, if $uv \in E$ then
there is a matching between all but one of the
neighbors of $u$ and $v$.)
Since $v$ has $2d-1-k$ edges to ${\mathcal O} \setminus (B \cup \{w\})$ and
$w$ has $2d-1-\ell$ edges to ${\mathcal E} \setminus (A \cup \{v\})$ we have
that $|C| = 4d -(k+\ell) \geq 2d$. From this
(\ref{inq-number.of.remaining.comps}) follows.

Inserting (\ref{inq-number.of.remaining.comps}) into
(\ref{inq-small.count}) and bounding $|\partial^\star A|$ and
$|\partial^\star B|$ by the maximum values of $|A|$ and $|B|$ (valid
since $T \subseteq {\mathcal E}$ (or ${\mathcal O}$) satisfies $|T| \leq
|\partial_{\rm ext}T|$, so $|\partial^\star T| \leq |T|$) and with the
remaining inequalities justified below, we have
\begin{eqnarray}
|{\mathcal C}_3^{\rm small}| & \leq & \exp_2\left\{\frac{n^d}{2}\left(\rho +
\frac{1}{d^{1/2}} + \frac{1}{d}\right)\right\} \cdot \sum_{i
\leq n^d/4d^{1/2}} {n^d/2 \choose i} \cdot
\sum_{j \leq (\rho + 1/2d^{1/2})n^d/2} {n^d/2 \choose j} \nonumber \\
& \leq & \exp_2\left\{\frac{n^d}{2}\left(\rho + \frac{1}{d^{1/2}} + \frac{1}{d} + H\left(\frac{1}{2d^{1/2}}\right) + H\left(\rho + \frac{1}{2d^{1/2}}\right)\right)\right\} \label{int4} \\
& \leq & \exp_2\left\{\frac{n^d}{2}\left(1-\Omega(1)\right)\right\} \label{int5}
\end{eqnarray}
for sufficiently large $d=d(\rho)$. In (\ref{int4}) we use the bound
$\sum_{i=0}^{[\beta M]}{M\choose i} \leq
2^{H(\beta)M}$ for $\beta \leq \frac{1}{2}$; in (\ref{int5}) we use
$H(\rho)+\rho < 1$. Using $2^{n^d/2} \leq |{\mathcal C}_3|$,
the lemma follows. \qed

\medskip

We now consider
$$
{\mathcal C}_3^{\rm large,~ even}:=({\mathcal C}_3^{b,\rho} \setminus {\mathcal C}_3^{\rm small}) \cap
{\mathcal C}_3^{\rm even}.
$$
By Lemma \ref{lem-bounding_small} and ${\mathcal E}$-${\mathcal O}$ symmetry, Theorem
\ref{thm-main-mixing} reduces to bounding (say)
\begin{equation} \label{inq-remaining}
\pi_3({\mathcal C}_3^{\rm large,~ even}) \leq
\exp\left\{-\frac{3n^{d-1}}{d^4\log^2 n}\right\}.
\end{equation}

Let ${\mathcal C}_3^{\rm large,~ even,~ nt}$ be the set of
$\chi \in {\mathcal C}_3^{\rm large,~ even}$ such that there is a
$\gamma \in \Gamma(\chi)$ with $|\gamma|\geq n^{d-1}$
and let
${\mathcal C}_3^{\rm large,~ even,~ triv} = {\mathcal C}_3^{\rm large,~ even} \setminus
{\mathcal C}_3^{\rm large,~ even,~ nt}$. We assert that
\begin{equation}\label{large2}
\pi_3({\mathcal C}_3^{\rm large,~ even,~ nt}) \leq
\exp\left\{-\Omega\left(\frac{n^{d-1}}{d}\right)\right\}
\end{equation}
and
\begin{equation}\label{large3}
\pi_3({\mathcal C}_3^{\rm large,~ even,~ triv})  \leq
\exp\left\{-\frac{4n^{d-1}}{d^4 \log^2 n}\right\};
\end{equation}
this gives (\ref{inq-remaining}) and so completes the proof of
Theorem \ref{thm-main-mixing}. Both (\ref{large2}) and (\ref{large3}) are
corollaries of Lemma \ref{lem-volume.bounds.from.gk}, and the steps
are identical to those that are used to bound the measures of
${\mathcal I}_{large, even}^{non-trivial}$ and ${\mathcal I}_{large, even}^{trivial}$
in \cite[Section 3.3]{Galvin}. We now give the details.

With the sum below running over all profiles $\ul{p}$ of the form
$(c,v)$ with $v \in V$ and $c \geq n^{d-1}$, and with the
inequalities justified below, we have
\begin{eqnarray}
\pi_3({\mathcal C}_3^{\rm large,~ even,~ nt})  & \leq &  \sum_{\ul{p}}
\pi_3({\mathcal C}_3(\ul{p})) \nonumber \\
& \leq &
n^{2d}\exp\left\{-\Omega\left(\frac{n^{d-1}}{d}\right)\right\} \label{large1} \\
& \leq & \exp\left\{-\Omega\left(\frac{n^{d-1}}{d}\right)\right\}, \nonumber
\end{eqnarray}
giving (\ref{large2}). We use Lemma \ref{lem-volume.bounds.from.gk}
in (\ref{large1}). The factor of $n^{2d}$ is for the choices of $c$
and $v$.

The verification of (\ref{large3}) involves finding an $i \in
[\Omega(\log d), O(d\log n)]$ and a set $\Gamma_i(\chi) \subseteq
\Gamma(\chi)$ of cutsets with the properties that $|\Gamma_i(\chi)|
\approx n^d/2^i$, $|\gamma| \approx 2^i$ for each $\gamma \in
\Gamma_i(\chi)$ and $\sum_{\gamma \in \Gamma_i(\chi)} |\gamma| \approx
n^{d-1}$. The measure of ${\mathcal C}_3^{\rm large,~ even,~ triv}$ is then at most
the product of a term that is exponentially small in $n^{d-1}$ (from
Lemma \ref{lem-volume.bounds.from.gk}), a term corresponding to the
choice of a fixed vertex in each of the interiors, and a term
corresponding to the choice of the collection of cutset sizes. The second
term will be negligible because $\Gamma_i(\chi)$ is small and the third
will be negligible because all $\gamma \in \Gamma_i(\chi)$ have similar
sizes.

More precisely, for $\chi \in {\mathcal C}_3^{\rm large,~ even,~ triv}$ and $\gamma \in
\Gamma(\chi)$ we have $|\gamma| \geq |{\rm int} \gamma|^{1-1/d}$ (by
(\ref{contour.prop.8}))
and
so
$$
\sum_{\gamma \in \Gamma(\chi)} |\gamma|^{d/(d-1)} \geq \sum_{\gamma \in
\Gamma(\chi)} |{\rm int} \gamma|
\geq |I^{\mathcal E}|
\geq  n^d/4d^{1/2}.
$$
The second inequality is from (\ref{cutset.conditions}) and the
third follows since $\chi \not \in {\mathcal C}_3^{\rm small}$.

Set $\Gamma_i(\chi) = \{\gamma \in \Gamma(\chi):2^{i-1} \leq |\gamma| <
2^i\}$. Note that $\Gamma_i(\chi)$ is empty for $2^i < d^2$ (again
by (\ref{contour.prop.8}))
and for $2^{i-1} > n^{d-1}$ so we may assume
that
\begin{equation} \label{inner.property.1}
2 \log d \leq i \leq (d-1)\log n+1.
\end{equation}
Since $\sum_{m=1}^\infty 1/m^2 = \pi^2/6$, there is an $i$ such that
\begin{equation} \label{quad.scale}
\sum_{\gamma \in \Gamma_i(\chi)} |\gamma|^{\frac{d}{d-1}} \geq
\Omega\left(\frac{n^d}{d^{1/2}i^2}\right).
\end{equation}
Choose the smallest such $i$ and set $\ell=|\Gamma_i(\chi)|$. We have
$\sum_{\gamma \in \Gamma_i(\chi)} |\gamma| \geq \Omega(\ell 2^i)$ (this
follows from the fact that each $\gamma \in \Gamma_i(\chi)$ satisfies
$|\gamma| \geq 2^{i-1}$) and
\begin{equation} \label{inner.property.2}
O\left(\frac{dn^d}{2^i}\right) \geq \ell \geq
\Omega\left(\frac{n^d}{2^{\frac{id}{d-1}}i^2d^{1/2}}\right).
\end{equation}
The first inequality follows from that fact that $\sum_\gamma |\gamma|
\leq dn^d=|E|$; the second follows from (\ref{quad.scale}) and the
fact that each $\gamma$ has $|\gamma|^{d/(d-1)} \leq 2^{di/(d-1)}$. We
therefore have $\chi \in {\mathcal C}_3(\ul{p})$ for some $\ul{p}=(c_1, v_1,
\ldots, c_\ell, v_\ell)$ with $\ell$ satisfying
(\ref{inner.property.2}), with
\begin{equation} \label{inner.property.3}
\sum_{j=1}^\ell c_j \geq O(\ell 2^i),
\end{equation}
with
\begin{equation} \label{inner.property.4}
c_j \leq 2^i
\end{equation}
for each $j$ and with $i$ satisfying (\ref{inner.property.1}). With
the sum below running over all $\ul{p}$ satisfying
(\ref{inner.property.1}), (\ref{inner.property.2}),
(\ref{inner.property.3}) and (\ref{inner.property.4}) we have
\begin{eqnarray}
\pi_3({\mathcal C}_3^{\rm large,~ even,~ triv}) & \leq & \sum_{\ul{p}}
\pi_3({\mathcal C}_3(\ul{p})). \label{largeeventriv}
\end{eqnarray}
The right-hand side of (\ref{largeeventriv}) is, by Lemma
\ref{lem-volume.bounds.from.gk}, at most
$$
d\log n~\max \left\{2^{\ell
i}{n^d \choose \ell}\exp\left\{-\Omega\left(\frac{\ell
2^i}{d}\right)\right\} : i~\mbox{satisfying (\ref{inner.property.1})}\right\}.
$$
The factor of $d\log n$ is an upper bound on the number of choices
for $i$; the factor of $2^{\ell i}$ is for the choice of the
$c_j$'s; and the factor ${n^d \choose \ell}$ is for the choice of
the $\ell$ (distinct) $v_j$'s. By (\ref{inner.property.1}) and the
second inequality in (\ref{inner.property.2}) we have (for $d$
sufficiently large)
$$
2^{\ell i}{n^d \choose \ell}  \leq   2^{\ell i}\left(\frac{n^d}{\ell}\right)^\ell
 \leq  2^{\ell i} \left(O\left(2^\frac{id}{d-1}i^2 d^{1/2}\right)\right)^\ell
 \leq  2^{4\ell i}
 = \exp\left\{o\left(\frac{2^i}{d}\right)\right\},
$$
so that in fact the right-hand side of (\ref{largeeventriv}) is at
most
$$
d \log n ~\max_i
\exp\left\{-\Omega\left(\frac{2^i\ell}{d}\right)\right\}.
$$
Taking $\ell$ as small as possible we see that this is at most
$$
d \log n ~\max_i \exp\left\{-\Omega\left(\frac{2^in^d}{d
2^{\frac{id}{d-1}}i^2d^{1/2}}\right)\right\}
$$
and taking $i$ as large as possible we see that it is at most
$\exp\{-4n^{d-1}/d^4 \log^2 n\}$. Putting these observation together
we obtain (\ref{large3}).

\section{Proof of Lemmas \ref{lem-volume.bounds.from.gk-pt} and \ref{lem-volume.bounds.from.gk}}
\label{sec-lemmaproofs}

In this section we complete the proofs of Theorems \ref{thm-influence1} and \ref{thm-glauber.slow.mixing} by establishing the two technical statements concerning cutsets from Sections \ref{sec-influenceproof} and \ref{sec-mixingproof}. Much of what follows is modified from \cite{Galvin} and \cite{GalvinKahn}. Because the cutsets described in Sections \ref{sec-cutsets} and \ref{sec-cutsets2} are quite similar, the two proofs proceed almost identically, and we give them in parallel. Before beginning this process we reduce Lemma \ref{lem-volume.bounds.from.gk} to (\ref{induction}) below. Let $\ul{p}=(c_0,v_0, \ldots, c_\ell,
v_\ell)$ be given. Set $\ul{p'}=(c_1 ,v_1, \ldots, c_\ell, v_\ell)$.
We will show
\begin{equation} \label{induction}
\frac{|{\mathcal C}_3(\ul{p})|}{|{\mathcal C}_3(\ul{p'})|} \leq
\exp\left\{-\Omega\left(\frac{c_0}{d}\right)\right\}
\end{equation}
from which Lemma \ref{lem-volume.bounds.from.gk} follows by a telescoping product. To obtain
(\ref{induction}) we define a one-to-many map $\varphi$ from
${\mathcal C}_3(\ul{p})$ to ${\mathcal C}_3(\ul{p'})$. We then define a flow
$\nu:{\mathcal C}_3(\ul{p}) \times {\mathcal C}_3(\ul{p'}) \rightarrow [0,\infty)$
supported on pairs $(\chi,\chi')$ with $\chi' \in \varphi(\chi)$
satisfying
\begin{equation}\label{eq-flow.out}
\forall \chi \in {\mathcal C}_3(\ul{p}),  \sum_{\chi' \in \varphi(\chi)}
\nu(\chi,\chi') =1
\end{equation}
and
\begin{equation}\label{eq-flow.in}
\forall \chi' \in {\mathcal C}_3(\ul{p'}), \sum_{\chi \in
\varphi^{-1}(\chi')} \nu(\chi,\chi') \leq
\exp\left\{-\Omega\left(\frac{c_0}{d}\right)\right\}.
\end{equation}
This easily gives (\ref{induction}).
To obtain Lemma \ref{lem-volume.bounds.from.gk-pt}, we prove a variant of (\ref{induction}) with ${\mathcal C}_3(\ul{p'})$ replaced by ${\mathcal C}_3^{\mathcal O}$ and ${\mathcal C}_3(\ul{p})$ replaced by ${\mathcal C}_3^{\mathcal O}(c_0,v_0)$.

In what follows, we write ${\mathcal D}$ for both ${\mathcal C}_3(\ul{p'})$ and ${\mathcal C}_3^{\mathcal O}$, and ${\mathcal C}$ for both ${\mathcal C}_3(\ul{p})$ and ${\mathcal C}_3^{\mathcal O}(c_0,v_0)$, and we use $V$ both for the vertex set of ${\mathbb Z}^d_n$ and that of $\Lambda$.

\medskip

For each $s \in \{\pm 1, \ldots, \pm d\}$, define $\sigma_s$, the
{\em shift in direction $s$}, by $\sigma_s(x)=x+e_s$, where $e_s$ is
the $s$th standard basis vector if $s>0$ and $e_s=-e_{-s}$ if $s<0$.
For $X \subseteq V$ write $\sigma_s(X)$ for $\{\sigma_s(x):x \in
X\}$. For $\gamma \in {\mathcal W}$ set $W^s = \{x \in
\partial_{\rm int}W:\sigma_{-s}(x) \not \in W\}$.

Let $\chi \in {\mathcal C}$ be given. For Lemma \ref{lem-volume.bounds.from.gk}, arbitrarily pick $\gamma \in
\Gamma(\chi) \cap {\mathcal W}(c_0,v_0)$ and set $W={\rm int} \gamma$. For Lemma \ref{lem-volume.bounds.from.gk-pt}, simply take $\gamma=\gamma(\chi)$ and $W=W(\gamma)$. Write $f$ for
the map from $\{0,1,2\}$ to $\{0,1,2\}$ that sends $0$ to $0$ and
transposes
$1$ and $2$. For each $s \in \{\pm 1, \ldots, \pm d\}$ and $S
\subseteq W$ define the function $\chi^s_S:V \rightarrow \{0,1,2\}$
by
$$
\chi^s_S(v) = \left\{ \begin{array}{ll}
                 0 & \mbox{if $v \in S$} \\
                 \chi(v) & \mbox{if $v \in (W^s\setminus S)
                 \cup (V \setminus W)$} \\
                 f(\chi(\sigma_{-s}(v))) & \mbox{if $v \in W \setminus
                 W^s$} \\
                      \end{array}
              \right.
$$
and set $\varphi_s(\chi)=\{\chi^s_S:S \subseteq W^s\}$.

\begin{claim} \label{claim-shift.works}
$\varphi_s(\chi) \subseteq {\mathcal D}$.
\end{claim}

\smallskip

\noindent {\em Proof}:
We begin with the observation that the graph $\partial_{\rm int} W \cup
\partial_{\rm ext} W$ is bipartite with bipartition
$(\partial_{\rm int} W, \partial_{\rm ext} W)$. This follows from
(\ref{contour.prop.2}). By (\ref{contour.prop.3}), $I \cap
(\partial_{\rm int} W \cup
\partial_{\rm ext} W)=\emptyset$ and so for each component $U$ of $\partial_{\rm int} W \cup
\partial_{\rm ext} W$, $\chi$ is identically $1$ on one of $U \cap \partial_{\rm int} W$, $U \cap \partial_{\rm ext} W$ and identically $2$ on the other.

Our main task is to show that $\varphi_s(\chi) \subseteq {\mathcal C}_3$; that is, that for any $S \subseteq W^s$ and edge $uv$, $\chi^s_S(u) \neq \chi^s_S(v)$. We consider several
cases.

If $u, v \not \in W$ then $\chi^s_S(u)=\chi(u)$ and $\chi^s_S(v) =
\chi(v)$. But $\chi(u) \neq \chi(v)$, so $\chi^s_S(u) \neq
\chi^s_S(v)$ in this case.

If $u \in W$ and $v \not \in W$ then $\chi^s_S(v) = \chi(v)$ and
$\chi^s_S(u) \in \{0,\chi(u)\}$ (we will justify this in a moment).
Since $v \in \partial_{\rm ext} W$ we have $\chi(v) \neq 0$ and we
cannot ever have $\chi(v)=\chi(u)$, so $\chi^s_S(u) \neq
\chi^s_S(v)$ in this case. To see that $\chi^s_S(u) \in
\{0,\chi(u)\}$, we consider subcases. If $u \in S$ then
$\chi^s_S(u)=0$. If $u \in W^s \setminus S$ then
$\chi^s_S(u)=\chi(u)$. Finally, if $u \in W \setminus W^s$ then
$\chi^s_S(u)=f(\chi(\sigma_{-s}(u)))$; and $f(\chi(\sigma_{-s}(u)))$ is
either $0$ or $\chi(u)$ depending on whether $\chi(\sigma_{-s}(u))$
equals $0$ or $\chi(v)$ ($\chi(\sigma_{-s}(u))$ cannot equal~$\chi(u)$).

If $u, v \in W\setminus W^s$ then $\chi^s_S(u)=f(\chi(\sigma_{-s}(u)))$
and $\chi^s_S(v)=f(\chi(\sigma_{-s}(v)))$. Since $f$ is a bijection and
$\chi(\sigma_{-s}(u)) \neq \chi(\sigma_{-s}(v))$ we have $\chi^s_S(u) \neq
\chi^s_S(v)$ in this case.

If $u \in W\setminus W^s$ and $v \in W^s \setminus S$ then
$\chi^s_S(u) \in \{0,\chi(u)\}$ (as in the second case above) and
$\chi^s_S(v)=\chi(v)$. Since $\chi(v) \neq 0$, we have $\chi^s_S(u)
\neq \chi^s_S(v)$.

Noting that it is not possible to have both $u, v \in W^s$, we
finally treat the case where $u \in W\setminus W^s$ and $v \in S$.
In this case $\chi^s_S(v)=\chi(v)=0$. Suppose (for a contradiction)
that $\chi^s_S(u)=0$. This can only happen if $\chi(\sigma_{-s}(u))=0$.
If $\sigma_{-s}(u)=v$, we have a contradiction immediately. Otherwise,
we have $\sigma_{-s}(v) \not \in W$ and so (since
$\sigma_{-s}(u)\sigma_{-s}(v) \in E$) $\sigma_{-s}(u) \in
\partial_{\rm int} W$, also a contradiction.

This verifies $\varphi_s(\chi) \subseteq {\mathcal C}_3$. We now verify that $\varphi_s(\chi) \subseteq {\mathcal D}$. In the setting of Lemma \ref{lem-volume.bounds.from.gk} this is true because $W$ is disjoint from the interiors of the remaining
cutsets in $\Gamma(\chi)$ and the operation that creates the elements
of $\varphi_s(\chi)$ only modifies $\chi$ inside $W$. In the setting of Lemma \ref{lem-volume.bounds.from.gk-pt} it follows from the fact that $W \cap \partial_{\rm int} \Lambda = \emptyset$. \qed

\begin{claim} \label{claim-unique.reconstruction}
Given $\chi' \in \varphi_s(\chi)$, $\chi$ can be uniquely
reconstructed from $W$ and $s$.
\end{claim}

\smallskip

\noindent {\em Proof}: We may reconstruct $\chi$ via
$$
\chi(v) = \left\{ \begin{array}{ll}
                 \chi'(v) & \mbox{if $v \in V \setminus W$} \\
                 f(\chi'(\sigma_s(v))) & \mbox{if $v \in W$}. \\
                      \end{array}
              \right.
$$
\qed

\medskip

We define the one-to-many map $\varphi$ from ${\mathcal C}$ to
${\mathcal D}$ by setting $\varphi(\chi)=\varphi_s(\chi)$ for a
particular direction $s$.
To define $\nu$ and $s$, we employ the notion of approximation also
used in \cite{GalvinKahn} and based on ideas introduced by
Sapozhenko in \cite{Sap}.
For $\gamma \in {\mathcal W}$, we say $A \subseteq V$ is an {\em approximation}
of $\gamma$ if
$$
A^{\mathcal E} \supseteq W^{\mathcal E} ~~~\mbox{and}~~~A^{\mathcal O} \subseteq W^{\mathcal O},
$$
$$
d_{A^{\mathcal O}}(x) \geq 2d-\sqrt{d} ~~\mbox{for all $x \in A^{\mathcal E}$}
$$
and
$$
d_{{\mathcal E} \setminus A^{\mathcal E}}(x) \geq 2d-\sqrt{d} ~~\mbox{for
all $y \in {\mathcal O} \setminus A^{\mathcal O}$},
$$
where $d_X(x)=|\partial x \cap X|$. Note that
by (\ref{contour.prop.5}),
$W(\gamma)$ is an approximation of $\gamma$.

Before stating our main approximation lemma, it will be convenient
to further refine our partition of cutsets. To this end set
$$
{\mathcal W}(w_e,w_o,v_0) = \left\{\gamma:
\mbox{$\gamma \in {\mathcal W}$ with $|W^{\mathcal O}|=w_o$, $|W^{\mathcal E}|=w_e$ and $v_0 \in W^{\mathcal E}$}
\right\}.
$$
Note that
by (\ref{contour.prop.2}) we have
$|\gamma| =2d(|W^{\mathcal O}|-|W^{\mathcal E}|)$ so ${\mathcal W}(w_e,w_o,v_0) \subseteq
{\mathcal W}((w_o-w_e)/2d, v_0)$.

\begin{lemma} \label{lem-k.comp.approx}
For each $w_e$, $w_o$ and $v_0$ there is a family ${\mathcal A}(w_e,w_o,v_0)$ of subsets of $V$
satisfying
$$
|{\mathcal A}(w_e,w_o,v_0)| \leq
\exp\left\{O\left((w_o-w_e)d^{-\frac{1}{2}}\log^{\frac{3}{2}}d\right)\right\}
$$
and a map $\pi:{\mathcal W}(w_e,w_o,v_0) \rightarrow {\mathcal A}(w_e,w_o,v_0)$ such that
for each $\gamma \in {\mathcal W}(w_e,w_o,v_0)$, $\pi(\gamma)$ is an approximation
of $\gamma$.
\end{lemma}

\smallskip

\noindent {\em Proof}: In the setting of Lemma \ref{lem-volume.bounds.from.gk-pt}, this is exactly \cite[Lemma 2.18]{GalvinKahn}; in the setting of Lemma \ref{lem-volume.bounds.from.gk} it is \cite[Lemma 4.2]{Galvin}.
\qed

\medskip

In both settings, the proof proceeds along the same lines. We begin by associating with each cutset a small set of vertices (much smaller than the size of the cutset) which weakly approximates the cutset in the sense that the neighborhood of the associated set separates the interior of the cutset from the exterior. This part of the proof combines algorithmic and probabilistic elements, and relies heavily on the structure of the lattice. The total number of weak approximations that can arise as we run over all cutsets of a given size is controlled in part by the fact that these weak approximations are connected (in a suitable sense); this property is inherited from the connectivity of the cutsets themselves. The second part of the proof proceeds by refining the weak approximations into approximations in the sense defined above. This part of the proof is purely algorithmic and uses no
properties of the lattice other than that it is regular and bipartite.

\medskip

We are now in a position to define $\nu$ and $s$. Recall that we have fixed, for each $\chi \in \varphi^{-1}(\chi')$, a particular cutset $\gamma$.
Our plan is to fix $w_e, w_o$ and $A \in
{\mathcal A}(w_e, w_o,v_0)$
and to consider the contribution to the sum in (\ref{eq-flow.in})
from those $\chi \in \varphi^{-1}(\chi')$ with $\pi(\gamma)=A$.
We will try to define $\nu$ in such a way that each of these
individual contributions to (\ref{eq-flow.in}) is small; to succeed
in this endeavor we must first choose $s$ with care. To this end,
given $A \in {\mathcal A}(w_e,w_o,v_0)$ set
$$
Q^{\mathcal E} = A^{\mathcal E} \cap \partial_{\rm ext}({\mathcal O} \setminus A^{\mathcal O})
~~\mbox{and}~~Q^{\mathcal O} = ({\mathcal O} \setminus A^{\mathcal O}) \cap
\partial_{\rm ext} A^{\mathcal E}.
$$

To motivate the introduction of $Q^{\mathcal E}$ and
$Q^{\mathcal O}$, note that for $\gamma \in \pi^{-1}(A)$ we have
(by (\ref{contour.prop.5}) and the definition of approximation)
$$
A^{\mathcal E}\setminus Q^{\mathcal E} \subseteq W^{\mathcal E},
$$
$$
{\mathcal E} \setminus A^{\mathcal E} \subseteq
{\mathcal E} \setminus W^{\mathcal E},
$$
$$
A^{\mathcal O}  \subseteq  W^{\mathcal O},
$$
and
$$
{\mathcal O} \setminus (A^{\mathcal O} \cup Q^{\mathcal O}) \subseteq  {\mathcal O} \setminus W^ {\mathcal O}.
$$
It follows that for each $\gamma \in \pi^{-1}(A)$, $Q^{\mathcal E} \cup Q^{\mathcal O}$
contains all
vertices whose location in the partition
$V = W \cup \overline{W}$ is as yet unknown.

We choose
$s(\chi)$ to be the smallest $s$ for which both of $|W^s| \geq
.8(w_o-w_e)$ and $|\sigma_s(Q^{\mathcal E}) \cap Q^{\mathcal O}| \leq 5|W^s|/\sqrt{d}$
hold.
This is the direction that minimizes the uncertainty
to be resolved when we attempt to reconstruct
$\chi$ from the partial information provided by $\chi' \in
\varphi^{-1}(\chi)$, $s$ and $A$.
(That such an $s$ exists is established in \cite[(49) and
(50)]{GalvinKahn} by an easy averaging argument).
Note that $s$ depends on $\gamma$ but not~$I$.

Now for each $\chi \in {\mathcal C}$ let $\gamma \in \Gamma(\chi)$ be a
particular cutset with $\gamma \in {\mathcal W}(c_0, v_0)$. Let $\varphi(\chi)$
be as defined before, with $s$ as specified above.
Define
$$
C = W^s \cap A^{\mathcal O} \cap \sigma_s(Q^{\mathcal E})
$$
and
$$
D = W^s \setminus C,
$$
and for each $\chi' \in \varphi(\chi)$ set
$$
\nu(\chi,\chi') = \left(\frac{1}{4}\right)^{|C \cap I(\chi')|}
\left(\frac{3}{4}\right)^{|C \setminus
I(\chi')|}\left(\frac{1}{2}\right)^{|D|}.
$$
Note that for $\chi \in \varphi^{-1}(\chi')$, $\nu(\chi,\chi')$
depends on $W$ but not on $\chi$ itself.

Since $C \cup D$ partitions $W$ we easily have (\ref{eq-flow.out}).
To obtain (\ref{induction}) we
must establish (\ref{eq-flow.in}).

Fix $w_e$, $w_o$ such that $2d(w_o-w_e)=c_0$. Fix $A \in
{\mathcal A}(w_e,w_o,v_0)$ and $s \in \{\pm 1, \ldots, \pm d\}$. For $\chi$
with $\gamma \in {\mathcal W}(w_e,w_o, v_0)$ write $\chi \sim_s A$ if it holds
that $\pi(\gamma)=A$ and $s(\chi)=s$. We claim that with $A, s, w_o$
and $w_e$ fixed, for $\chi' \in {\mathcal D}$
\begin{equation} \label{inq-black.box}
\sum \left\{\nu(\chi,\chi'):\chi \sim_s A,~\chi \in
\varphi^{-1}(\chi')\right\} \leq
\left(\frac{\sqrt{3}}{2}\right)^{w_o-w_e}.
\end{equation}

We now describe the proof of (\ref{inq-black.box}).
Write ${\mathcal C}(w_e,w_o,s,A,\chi')$ for the set of all $\chi
\in {\mathcal C}$ such that $W \in {\mathcal W}(w_e,w_o,v_0)$,
$\pi(\gamma)=A$, $s(\chi)=s$ and $\chi' \in \varphi(\chi)$ and set $U=
Q^{\mathcal E} \cap \sigma_{-s}(\chi')$. Say that a triple $(K,L,M)$ is {\em
good} for $\chi$ if it satisfies the following conditions.
$$
\mbox{$K \cup L \cup M$ is a minimal vertex cover of $Q^{\mathcal E} \cup
Q^{\mathcal O}$,}
$$
$$
\mbox{$K \subseteq Q^{\mathcal O}$, $L \subseteq U$ and $M \subseteq Q^{\mathcal E}
\setminus U$}
$$
and
$$
\mbox{$K=\partial_{\rm ext}(U \setminus L)$}.
$$
We begin by establishing that $\chi \in
{\mathcal C}(w_e,w_o,s,A,\chi')$ always has a good triple.

\begin{lemma} \label{lem-core-exists}
For each $\chi \in {\mathcal C}(w_e,w_o,s,A,\chi')$ the triple
$$
(\hat{K}, \hat{L}, \hat{M}):=(W \cap Q^{\mathcal O}, U \setminus W, (Q^{\mathcal E}
\setminus U) \setminus W)
$$
is good for $\chi$.
\end{lemma}

\smallskip

\noindent {\em Proof}: \cite[around discussion of (54)]{GalvinKahn}.
\qed

\medskip

In view of Lemma \ref{lem-core-exists} there is a triple $(K, L, M)$
that is good for $\chi$ and which has $|K|+|L|$ as small as
possible. Choose one such, say $(K_0(\chi), L_0(\chi), M_0(\chi))$.
Set $K'(\chi)=K_0\setminus \hat{K}$ and $L'(\chi)=L_0\setminus
\hat{L}$. Lemma \ref{lem-unlabeled_main_bound} below establishes an
upper bound on $\nu(\chi,\chi')$ in terms of $|K_0|$, $|L_0|$,
$|K'|$ and $|L'|$, and Lemma \ref{lem-reconstruction} shows that for
each choice of $K'$, $L'$ there is at most one $\chi$ contributing
to the sum in the lemma. These two lemmas combine to give
(\ref{inq-black.box}).
\begin{lemma} \label{lem-unlabeled_main_bound}
For each $\chi \in {\mathcal C}(w_e,w_o,s,A,\chi')$,
\begin{eqnarray*}
\nu(\chi,\chi') & \leq & \left(\frac{\sqrt{3}}{2}\right)^{w_o-w_e}
\frac{2^{|K_0|}}{3^{|K_0|+|L_0|}
2^{|K'|-|L'|}} \\
& := & B(K',L').
\end{eqnarray*}
\end{lemma}

\smallskip

\noindent {\em Proof}: We follow \cite[from just before (55) to just
after (60)]{GalvinKahn}, making superficial changes of notation.
\qed

\medskip

The inequality in Lemma \ref{lem-unlabeled_main_bound} is the
$3$-coloring analogue of the main inequality of \cite{GalvinKahn}.
The key observation that makes this inequality useful is the
following.
\begin{lemma} \label{lem-reconstruction}
For each $w_e$, $w_o$, $s$, $A$, $\chi'$, $K'$ and $L'$, there is at
most one $\chi$ with $\chi \in {\mathcal C}(w_e,w_o,s,A,\chi')$,
$K'=K'(\chi)$ and $L'=L'(\chi)$.
\end{lemma}

\smallskip

\noindent {\em Proof}: In \cite[(56) and following]{GalvinKahn} it
is shown that $K'$ and $L'$ determine $W^{\mathcal O}$ via
$$
\hat{K}= (K_0\setminus K') \cup (\partial_{\rm ext} L' \cap Q^{\mathcal O})
$$
and so $W$ (via $W^{\mathcal E}=\{v \in {\mathcal E}:\partial v \subseteq W^{\mathcal O}\}$).
But then by Claim \ref{claim-unique.reconstruction} $K'$ and $L'$
determine $\chi$. \qed

\medskip

Lemmas \ref{lem-unlabeled_main_bound} and \ref{lem-reconstruction}
together now easily give (\ref{inq-black.box}):
\begin{eqnarray*}
\sum_{\chi \in {\mathcal C}(w_e,w_o,s,A,\chi')} \!\!\!\!\!
\nu(\chi,\chi') & \leq & \sum_{K' \subseteq K_0, ~L' \subseteq L_0}
\!\!\!B(K', L')
 \\
& \leq & \left(\frac{\sqrt{3}}{2}\right)^{w_o-w_e}.
\end{eqnarray*}

We have now almost reached (\ref{eq-flow.in}).
With the steps justified below we have that for each $\chi' \in
{\mathcal D}$
\begin{eqnarray}
\sum_{\chi \in \varphi^{-1}(\chi')}  \nu (\chi,\chi')  & \leq &
\sum \left\{\nu(\chi,\chi'): \chi \sim_s A, ~\chi \in
\varphi^{-1}(\chi')\right\} \nonumber \\
& \leq & 2d c_0^{\frac{2d}{d-1}}|{\mathcal A}(w_e, w_o,
v_0)|\left(\frac{\sqrt{3}}{2}\right)^{\frac{c_0}{2d}}
\label{overcount} \\
& \leq & 2d
c_0^{\frac{2d}{d-1}}\exp\left\{-\Omega\left(c_0/d\right)\right\}
\label{overcount2} \\
& \leq & \exp\left\{-\Omega\left(c_0/d\right)\right\},
\label{overcount3}
\end{eqnarray}
completing the proof of (\ref{eq-flow.in}). In the first inequality,
the sum on the right-hand side is over all choices of $w_e$, $w_o$, $s$ and  $A$. In
(\ref{overcount}), we note that there are $|{\mathcal A}(w_e, w_o, v_0)|$
choices for
$A$, $2d$ choices for $s$ and
$c_0^{d/(d-1)}$ choices for each of $w_e$, $w_o$ (this is because
$c_0\geq (w_e+w_o)^{1-1/d}$, by (\ref{contour.prop.8})), and we
apply (\ref{inq-black.box})
to bound the summand. In (\ref{overcount2}) we use Lemma
\ref{lem-k.comp.approx}. Finally in (\ref{overcount3}) we use $c_0
\geq d^2$ (again by (\ref{contour.prop.8})) to bound
$2dc_0^{2d/(d-1)}=\exp\{o(c_0 /d)\}$.

\section{Proof of Theorem \ref{thm-maxent} (measures of maximal entropy)} \label{sec-entropyproof}

Here we establish that the Gibbs measure studied in Theorem \ref{thm-influence1} is a measure of maximal entropy.
Recall that for a probability distribution $X$ with finite range that takes on value $x$ with probability $p(x)$, the {\em entropy} of $X$ is
$$
H(X) = - \sum_{x \in {\rm range}(X)} p(x)\log p(x).
$$
We have $H(X) \leq \log |{\rm range}(X)|$ with equality if and only if $X$ is uniform.

Let $\Lambda_n$ be the box $\{-n, \ldots, n\}^d$, and let ${\mathcal C}'_3(\Lambda_n)$ be the set of colorings of $\Lambda_n$ that can be extended to a coloring of ${\mathbb Z}^d$. The {\em topological entropy} of ${\mathcal C}_3$ (the set of $3$-colorings of ${\mathbb Z}^d$) is
$$
{\mathcal H}^{\rm topo}({\mathcal C}_3) = \lim_{n \rightarrow \infty} \frac{\log |{\mathcal C}'_3(\Lambda_n)|}{|\Lambda_n|}.
$$
Let $\mu$ be any measure on $({\mathcal C}_3,{\mathcal F}_{\rm cyl})$ and let $X_n$ be the restriction to $\Lambda_n$ of an element of ${\mathcal C}_3$ chosen according to $\mu$ (so the range of $X_n$ is a subset of ${\mathcal C}'_3$). The {\em measure-theoretic entropy} of ${\mathcal C}_3$ with respect to $\mu$ is
$$
{\mathcal H}^{\mu}({\mathcal C}_3) = \lim_{n \rightarrow \infty} \frac{H(X_n)}{|\Lambda_n|}.
$$
Note that ${\mathcal H}^{\mu}({\mathcal C}_3)$ is always at most ${\mathcal
H}^{\rm topo}({\mathcal C}_3)$. We say that $\mu$ is a {\em measure of maximal
entropy} if ${\mathcal H}^{\mu}({\mathcal C}_3) = {\mathcal H}^{\rm
topo}({\mathcal C}_3)$. The sense of measure of maximal entropy is that the
restriction of $\mu$ to any finite subset of ${\mathbb Z}^d$ is supported
(asymptotically) on as large a set as possible. (See e.g.\ \cite{BurtonSteif} for a more thorough discussion of these topics.)

We wish to show that $\mu^{\chi(0,{\mathcal O})}$ (as described in the introduction) is a measure of maximal entropy. Fix $m$ and $n$ satisfying $m > n$. Let $\mu_m=\mu_m^{\chi(0,{\mathcal O})}$ be as described in the introduction, and let
$X^m_n$ be the restriction to $\Lambda_n$ of a coloring chosen according to $\mu_m$.
We will show that
\begin{equation} \label{max-ent:to-show-1}
H(X^m_n) \geq \log |{\mathcal C}'_3(\Lambda_n)| -2|\partial_{\rm int} \Lambda_n|\log 3.
\end{equation}
This is enough to show that $\mu^{\chi(0,{\mathcal O})}$ is a measure of maximal entropy, since $|\partial_{\rm int} \Lambda_n|=o(\log |{\mathcal C}'_3(\Lambda_n)|)$.

Since for any random variable $X$ we have $H(X)\geq -\log \max_x p(x)$, we will have (31) if we show that, for each $\tau \in {\mathcal C}'_3(\Lambda_n)$, we have
\begin{equation}\label{max-ent:to-show-2}
\Pr(X^m_n=\tau) \leq \frac{3^{2|\partial_{\rm int} \Lambda_n|}}{|{\mathcal C}'_3(\Lambda_n)|}.
\end{equation}

We need the following lemma. Here $\Sigma$ is an arbitrary finite bipartite graph with bipartition ${\mathcal E}\cup {\mathcal O}$.
\begin{lemma}\label{Lcol}
Fix ${\mathcal E}'\subseteq {\mathcal E}$ and ${\mathcal O}'\subseteq{\mathcal O}$ arbitrarily and
let $\mu$ be uniform measure on ${\mathcal C}_3(\Sigma)$.
For any
${\mathcal E}''\subseteq {\mathcal E}\setminus {\mathcal E}'$, ${\mathcal O}''\subseteq{\mathcal O}\setminus {\mathcal O}'$,
$$
\mu(\mbox{$\chi \equiv 0$ on ${\mathcal E}''$ and $\chi \equiv 1$ on ${\mathcal O}''$}~|~\mbox{$\chi \equiv 0$ on ${\mathcal E}'$ and $\chi \equiv 1$ on ${\mathcal O}'$})
\geq 3^{-|{\mathcal E}'' \cup {\mathcal O}''|}.
$$
\end{lemma}

\smallskip

\noindent {\em Proof}: We proceed by induction on $|{\mathcal E}'' \cup {\mathcal O}''|$, beginning with the case $|{\mathcal E}'' \cup {\mathcal O}''|=1$. Without loss of generality, we may take ${\mathcal O}''=\emptyset$ and ${\mathcal E}''=\{x\}$ for some $x \in {\mathcal E} \setminus {\mathcal E}'$.
Write ${\mathcal C}'$ for the set of those $\chi$ satisfying $\chi|_{{\mathcal E}'} \equiv 0$ and $\chi|_{{\mathcal O}'}\equiv 1$, and, for $i \in \{0,1,2\}$, write ${\mathcal C}'_i$ for
$\{\chi\in {\mathcal C}':\chi(x)=i\}$. We wish to show that $|{\mathcal C}'_0|/|{\mathcal C}'|\geq 1/3$, for which (by $1$-$2$ symmetry) it is enough to show $|{\mathcal C}'_1| \leq |{\mathcal C}'_0|$.

To verify this last inequality, consider the following map from ${\mathcal C}'_1$ to ${\mathcal C}'$: for $\chi\in {\mathcal C}'_1$, let $C$ be the set of vertices in $\Sigma$ reachable from $x$ using only vertices
colored $0$ and $1$, and
let $\chi'$ be obtained from $\chi$ by interchanging $0$ and $1$
on $C$. We must have $C\cap ({\mathcal E}'\cup{\mathcal O}')=\emptyset$
(since otherwise we would have an odd path from $x$ to ${\mathcal E}'$ or
an even path from $x$ to ${\mathcal O}'$), so that in fact $\chi' \in {\mathcal C}'_0$.
Moreover, the map is injective
since we can recover $\chi$ by
interchanging $0$ and $1$ on the set of vertices in $\Sigma$ reachable from $x$ using only vertices
colored $0$ and $1$ (under $\chi'$).

For the induction step, consider the case $|{\mathcal E}'' \cup {\mathcal O}''|=t>1$ where without loss of generality $|{\mathcal E}''|>0$. Fix $x \in {\mathcal E}''$. We have
\begin{eqnarray*}
& \mu(\mbox{$\chi \equiv 0$ on ${\mathcal E}''$ and $\chi \equiv 1$ on ${\mathcal O}''$}~|~\mbox{$\chi \equiv 0$ on ${\mathcal E}'$ and $\chi \equiv 1$ on ${\mathcal O}'$}) = & \\
& \mu(\mbox{$\chi \equiv 0$ on ${\mathcal E}'' \setminus \{x\}$ and $\chi \equiv 1$ on ${\mathcal O}''$}~|~\mbox{$\chi \equiv 0$ on ${\mathcal E}'$ and $\chi \equiv 1$ on ${\mathcal O}'$})  ~ \times & \\
& \mu(\mbox{$\chi(x)=0$}~|~\mbox{$\chi \equiv 0$ on ${\mathcal E}' \cup ({\mathcal E}'' \setminus \{x\})$ and $\chi \equiv 1$ on ${\mathcal O}' \cup {\mathcal O}''$}).
\end{eqnarray*}
The first term in the product above is at least $1/3$ (it is another instance of the base case), and the second term is at least $3^{-(t-1)}$ (by induction), so the product is at least $3^{-t}$.
\qed

\medskip

Now let $A=W_m\setminus (\Lambda_n\setminus \partial_{\rm int}\Lambda_n)$ (recall from Section \ref{sec-intro} that $W_m$ is the box $\{-m, \ldots, m\}^d$ together with all of the odd vertices of the box $\{-(m+1), \ldots, m+1\}^d$).
For $\tau\in {\mathcal C}'_3(\Lambda_n)$, let
$N(\tau)$
be the number of $\chi \in {\mathcal C}_3(A)$
that agree with $\tau$ on $\Lambda_n$ and can be extended
to colorings in ${\rm supp}(\mu_m) := \{\chi\in {\mathcal C}_3:\mu_m(\chi)>0\}$.
Thus $N(\tau)$ depends only on the restriction of $\tau$ to
$\partial_{\rm int}\Lambda_n$, and $N(\tau)=\Pr(X^m_n=\tau)|{\rm supp}(\mu_m)|$.

Set
$$
{\mathcal C}_*(\Lambda_n) =
\{\tau_0\in {\mathcal C}'_3(\Lambda_n):\mbox{$\tau_0\equiv 0$ on $(\partial_{\rm int}\Lambda_n)\cap {\mathcal O}~$ and $~\tau_0\equiv 1$ on $(\partial_{\rm int} \Lambda_n)\cap {\mathcal E}$}\}.
$$
By Lemma \ref{Lcol} (with $\Sigma=A \cup \partial_{\rm ext} W_m$, ${\mathcal E}' \cup {\mathcal O}'= \partial_{\rm ext} W_m$ and ${\mathcal E}'' \cup {\mathcal O}''= \partial_{\rm int} \Lambda_n$) we have, for any $\tau_0\in {\mathcal C}_*(\Lambda_n)$ and $\tau\in {\mathcal C}'_3(\Lambda_n)$,
\begin{equation} \label{nurm}
N(\tau_0) \geq
3^{-|\partial_{\rm int} \Lambda_n|}N(\tau) = 3^{-|\partial_{\rm int}\Lambda_n|} \Pr(X^m_n=\tau)|{\rm supp}(\mu_m)|.
\end{equation}
Another application of Lemma \ref{Lcol} (with $\Sigma=\Lambda_n$, ${\mathcal E}' \cup {\mathcal O}'= \emptyset$ and ${\mathcal E}'' \cup {\mathcal O}''= \partial_{\rm int} \Lambda_n$) yields
$$
\frac{|{\mathcal C}_*(\Lambda_n)|}{|{\mathcal C}'_3(\Lambda_n)|} \geq 3^{-|\partial_{\rm int}\Lambda_n|},
$$
and so, since $N(\tau_0)/|{\rm supp}(\mu_m)| = \Pr(X^m_n=\tau_0) \leq |{\mathcal C}_*(\Lambda_n)|^{-1}$
we get
\begin{equation} \label{nurm2}
N(\tau_0) \leq \frac{3^{|\partial_{\rm int}\Lambda_n|}}{|{\mathcal C}'_3(\Lambda_n)|} |{\rm supp}(\mu_m)|.
\end{equation}
Combining (\ref{nurm}) and (\ref{nurm2}) we get (\ref{max-ent:to-show-2}).


\begin{thebibliography}{99}

\bibitem{ammb}
D. Achlioptas, M. Molloy, C. Moore, and F. Van Bussel, Sampling grid colourings with
fewer colours, {\em
Proc. LATIN '04},  80--89.

\bibitem{BollobasLeader2}
B. Bollob\'as and I. Leader, Edge-isoperimetric
inequalities in the grid, {\em Combinatorica} {\bf 11} (1991),
299--314.

\bibitem{BorgsChayesFriezeKimTetaliVigodaVu}
C. Borgs, J. Chayes, A. Frieze, J.H. Kim, P. Tetali, E. Vigoda, V.
Vu, Torpid Mixing of some Monte Carlo Markov Chain algorithms in
Statistical Physics, {\em Proc. IEEE FOCS '99}, 218--229.

\bibitem{BrightwellWinkler}
G. Brightwell and P. Winkler, Graph homomorphisms and phase
transitions, {\em J. Combin. Theory Ser. B} {\bf 77} (1999),
221--262.

\bibitem{BrightwellWinkler4}
G. Brightwell and P. Winkler,
Gibbs measures and dismantlable graphs, {\em J. Combin. Theory Ser. B} {\bf 78} (2000), 141--166.

\bibitem{bdg}
R. Bubley, M. Dyer and C. Greenhill, Beating the $2\Delta$ bound for
approximately counting colourings: a computer-assisted proof of
rapid mixing, {\em Proc. ACM-SIAM SODA '98}, 355--363.

\bibitem{BurtonSteif}
R. Burton and J. Steif,
Non-uniqueness of measures of maximal entropy for subshifts of finite type, {\em Ergodic Theory and Dynamical Systems} {\bf 14}
(1994), 213--235.

\bibitem{Dobrushin4}
R. Dobrushin, The description of a random field by means of conditional probabilities and conditions of its regularity, {\em Theory Probab. Appl.} {\bf 13} (1968), 197--224.

\bibitem{DyerFriezeJerrum}
M. Dyer, A. Frieze and M. Jerrum, On counting independent sets in
sparse graphs, {\em SIAM J. Comp.} {\bf 31} (2002), 1527--1541.

\bibitem{EngbersGalvin-tori}
J. Engbers and D. Galvin, $H$-coloring tori, {\em J. Combin. Theory Ser. B} {\bf 102} (2012), 1110--1133.
% DJG Oct 9 2012 - updated this reference

%\bibitem{fs}
%S.J. Ferreira and A.D. Sokal, Antiferromagnetic Potts model on the square lattice: a high precision
%Monte Carlo study, {\em J. Stat. Phys.} {\bf 96} (1999), 461--530.

\bibitem{FriezeVigoda}
A. Frieze and E. Vigoda, A survey on the use of Markov chains to
randomly sample colourings, in {\em Combinatorics,
Complexity and Chance}, Oxford University Press, 2007.

\bibitem{Galvin2}
D. Galvin, Sampling $3$-colourings of regular bipartite graphs, {\em Electron. J. Probab.} {\bf 12} (2007), 481--497.

\bibitem{Galvin}
D. Galvin, Sampling independent sets on the discrete torus, {\em Random Structures Algorithms} {\bf 33} (2008), 356--376.

\bibitem{GalvinKahn}
D. Galvin and J. Kahn, On phase transition in the hard-core model on
${\mathbb Z}^d$, {\em Comb. Prob. Comp.} {\bf 13} (2004), 137--164.

\bibitem{GalvinRandall}
D. Galvin and D. Randall, Torpid Mixing of Local Markov Chains on $3$-Colorings of the Discrete Torus, {\em Proc. ACM--SIAM SODA} (2007), 376--384.

\bibitem{GareyJohnson}
M. Garey and D. Johnson, {\em Computers and Intractability: A Guide to the Theory of NP-Completeness}, W. H. Freeman, San Francisco, 1979.

\bibitem{Georgii2}
H.-O. Georgii,
{\em Gibbs Measures and Phase Transitions},
de Gruyter, Berlin, 1988.

\bibitem{gmp}
L.A. Goldberg, R. Martin and M. Paterson,
Random sampling of 3-colourings in ${\mathbb Z}^2$,
{\em Random Structures Algorithms} {\bf 24} (2004), 279--302.

\bibitem{hv}
T. Hayes and E. Vigoda, Coupling with the stationary distribution and
improved sampling for colorings and independent sets,
{\em Proc. ACM-SIAM SODA '05},  971-979.

\bibitem{jer}
M.R. Jerrum,
A very simple algorithm for estimating the number of $k$-colorings
of a low-degree graph,
{\em Random Structures Algorithms} {\bf 7} (1995), 157--165.

\bibitem{JerrumSinclair}
M. Jerrum and A. Sinclair, The Monte Carlo Markov chain method: an
approach to approximate counting and integration, in {\em
Approximation Alorithms for NP-hard problems}, PWS, 1996.

\bibitem{Kotecky}
R. Koteck\'y, Long-range order for antiferromagnetic Potts models, {\em Phys. Rev. B} {\bf 31} (1985), 3088--3092.

\bibitem{RK}
R. Koteck\'y, personal communication.

\bibitem{LawlerSokal}
G.F. Lawler and A.D. Sokal,
Bounds on the $L_2$ spectrum for Markov chains and Markov processes: a generalization
of Cheeger's inequality,
{\em Trans. Amer. Math. Soc.} {\bf 309} (1988), 557-580.

\bibitem{lrs}
M.~Luby, D.~Randall, and A.J.~Sinclair,
Markov Chains for Planar Lattice Structures,
{\em SIAM J. Comput.} {\bf 31} (2001), 167--192.

\bibitem{Peierls}
R. Peierls, Ising's Model of Ferromagnetism, {\em Proc. Cambridge Philos. Soc.} {\bf 32} (1936), 477--481.

\bibitem{Peled}
R. Peled, High-Dimensional Lipschitz Functions are Typically Flat, arXiv:1005.4636, May 25 2010.

\bibitem{RP}
R. Peled, personal communication.


\bibitem{Potts}
R. Potts, Some generalized order-disorder transformations, {\em Proc. Cambridge Philos. Soc.} {\bf 48} (1952), 106--109.

\bibitem{RandallTetali}
D. Randall and P. Tetali, Analyzing Glauber dynamics by comparison of Markov chains, {\em Journal of Mathematical Physics} {\bf 41} (2000), 1598--1615.

\bibitem{Sap}
A.A. Sapozhenko, On the number of connected subsets with given
cardinality of the boundary in bipartite graphs, {\em Metody
Diskret. Analiz.} {\bf 45} (1987), 42--70.  (Russian.)

\bibitem{sinc}
A.J. Sinclair, {\em Algorithms for random generation \& counting: a
Markov chain approach,} Birkh\"{a}user, Boston, 1993.

\bibitem{Thomas}
L.E. Thomas,
Bound on the mass gap for finite volume stochastic Ising models at low temperature,
{\em Commun. Math. Phys.} {\bf 126} (1989), 1-11.

\bibitem{swk0}
J.S. Wang, R.H. Swendsen, and R. Koteck\'y, Antiferromagnetic Potts models, {\em Phys. Rev.
Lett.} {\bf 63} (1989) 109--112.

\bibitem{swk}
J.S. Wang, R.H. Swendsen, and R. Koteck\'y, Three-state
antiferromagnetic Potts models: A Monte Carlo study, {\em Phys. Rev.
B} {\bf 42} (1990) 2465--2474.

\end{thebibliography}
\end{document}